\let\ORIlabel\label
\let\ORIrefstepcounter\refstepcounter
\AddToHook{package/hyperref/before}{%
  \let\label\ORIlabel
  \let\refstepcounter\ORIrefstepcounter
}

\documentclass[final,reqno]{siamart190516}
\usepackage[utf8]{inputenc}

\title{Adaptive precision incomplete matrix factorization\thanks{Version of \today.}}

\newcommand{\readingaffil}{Department of Mathematics and Statistics, University of Reading, Whiteknights,
Reading RG6 6AX, UK and Scientific Computing Department, STFC Rutherford Appleton Laboratory, Didcot, Oxfordshire OX11 0QX, UK}

\ifdefined\LaTeXML
  \newcommand{\andrewthanks}{\thanks{\readingaffil\ (\email{andrew.lister@stfc.ac.uk})}}
  \newcommand{\jenniferthanks}{\thanks{\readingaffil\ (\email{jennifer.scott@stfc.ac.uk})}}
\else
  \newcommand{\andrewthanks}{\thanks{\readingaffil\ (\email{andrew.lister@stfc.ac.uk, jennifer.scott@stfc.ac.uk})}}
  \newcommand{\jenniferthanks}{\footnotemark[2]}
\fi

\author{Andrew Lister\andrewthanks
\and
Theo Mary%
  \thanks{Sorbonne Universit\'e, CNRS, LIP6, Paris, France
  (\email{theo.mary@lip6.fr})}
\and
Jennifer Scott\jenniferthanks
}

\headers{Adaptive precision incomplete matrix factorization}{Lister, Mary, and Scott}

\usepackage{amssymb}
\usepackage{amsmath}
\usepackage{algorithm}
\usepackage{algorithmic}
\usepackage{color}
\usepackage{soul}
\usepackage{adjustbox}
\usepackage{graphicx}
\usepackage{subcaption}
\usepackage{booktabs}
\usepackage{multirow}

\newsiamremark{remark}{Remark}

\definecolor{hotpink}{rgb}{0.9,0,0.5}
\definecolor{mygre}{RGB}{119,172,48}
\hypersetup{colorlinks,linkcolor=blue,citecolor=hotpink,urlcolor=blue}

\numberwithin{equation}{section}
\def\R{\mathbb{R}}
\def\norm#1{\|#1\|}
\def\norminf#1{\|#1\|}   

\def\abs#1{|#1|}
\def\epsabs{\varepsilon_a}
\def\epsrel{\varepsilon}
\def\ceps{c_\varepsilon}
\def\cepsp{c_{\varepsilon,p}}
\def\uh{\widehat{u}}
\def\ut{\widehat{u}'}
\def\lh{\widehat{\ell}}

\def\Uh{\widehat{U}}
\def\Lh{\widehat{L}}
\def\Rh{\widehat{R}}
\def\xh{\widehat{x}}
\def\yh{\widehat{y}}

\DeclareMathOperator{\fl}{\operatorname{f\kern.2ptl}} 
\DeclareMathOperator{\diag}{\operatorname{diag}}
\DeclareMathOperator{\op}{\operatorname{op}}

\usepackage{comment}

\usepackage{xspace}
\newcommand{\APILU}{AP-ILU\xspace}
\newcommand{\APIC}{AP-IC\xspace}
\begin{document}
\maketitle

\begin{abstract}
Incomplete factorizations are widely used to precondition Krylov subspace methods for solving large sparse linear systems.
Recent work has explored computing approximate factorizations in low precision, but using only one uniformly low precision limits the achievable accuracy and can degrade the quality of the preconditioner.
We propose adaptive precision incomplete factorization algorithms in which we store each entry of the factors in a precision chosen according to its magnitude and a target accuracy.
We present a rounding error analysis showing that the adaptive precision threshold-based incomplete LU and incomplete Cholesky factorizations achieve a backward error of the order of the target accuracy. We also discuss how to incorporate pivoting and scaling and how to perform the triangular solves with the adaptive precision factors.
Finally, we explain how the approach applies to memory-limited incomplete Cholesky factorizations by using the magnitude of the largest dropped entry as the target accuracy.
Numerical experiments on a wide range of matrices from the SuiteSparse collection show that adaptive precision can substantially reduce the storage required by the factors without degrading the quality of the preconditioner and improve the storage-weighted iteration count and overall storage required compared with uniform 64-, 32-, and 16-bit factorizations.
\end{abstract}

\begin{keywords}
  sparse matrices, incomplete factorizations, preconditioning, linear systems, floating-point arithmetic,
  mixed precision, adaptive precision, rounding error analysis
\end{keywords}

\begin{AMS}
65F05, 65F08, 65F10, 65F35, 65F50, 65G50, 65Y20  
\end{AMS}

\section{Introduction}

Our focus is on solving sparse linear systems of equations $Ax = b$, where $A \in \mathbb{R}^{n \times n}$ 
is nonsingular and $x, b \in \mathbb{R}^{n}$. In large-scale practical applications, $A$ can be too large for a direct solver to be feasible.
In such cases, a preconditioned Krylov subspace method such as GMRES or, for symmetric positive definite problems, the conjugate gradient (CG) method, is commonly employed. For a given linear system, determining an appropriate preconditioner can be challenging.
Incomplete matrix factorizations are among the most widely used algebraic preconditioners. The key idea is to compute an approximate LU or Cholesky factorization of the matrix $A$, discarding some entries that would occur in a complete factorization, thereby trading exactness for computational efficiency in terms of both time and storage. The resulting incomplete factors can be used to accelerate the Krylov subspace method by providing an approximation to the inverse. For an introduction to incomplete factorization
preconditioners, see   \cite{benz:02,saad:2003,sctu:2023} and the extensive reference lists therein.

In recent years, the emergence of low precision floating-point
arithmetic in computer hardware has led to a resurgence of interest in the use of mixed precision
numerical linear algebra; see, for example, the survey articles
\cite{abdel:2021,hima22,aabb}. This has included work on computing
incomplete factorizations~\cite{sctu:2024,sctu:2025a} or, more broadly, approximate factorizations~\cite{abhl23,abhl24,cakh23,hima19,hipr:2021}
in low precision and then employing them to precondition higher precision iterative methods such as iterative refinement, GMRES, or CG~\cite{hipr:2021,sctu:2024,sctu:2025a}. However, uniformly reducing
the precision of the factorization can impact the quality of the preconditioner and the convergence of the solver if the precision is too low,
especially for ill-conditioned problems, thereby limiting the use of lower precisions and the potential performance gains.

In this article, we propose \textit{adaptive precision} incomplete factorization algorithms,
which perform the operations within the factorization in the chosen compute precision and store different entries of the factors in potentially different floating-point precisions.
This adaptive precision paradigm~\cite[sect.~14]{hima22} has been successfully used for several linear algebra operations, such as low-rank matrix and tensor approximations~\cite{abbg23,bmp25,bmv}, sparse matrix--vector product~\cite{gjmm24,gjmm24a}, relaxed GMRES~\cite{agjm}, tile factorizations~\cite{amlk,cmal}, and 
block Jacobi~\cite{adfh19,facq21}, additive Schwarz~\cite{cjmn}, and approximate inverse~\cite{khca25} preconditioners, but has not, to our knowledge, been explored for incomplete factorizations.
We carry out an error analysis that shows that the criterion to choose the precision for storing an entry of the factors should be based
on both the magnitude of the entries and the target accuracy of the incomplete factorization.
Storing an entry in low precision introduces an error that is proportional to its magnitude,
so small entries can be stored in correspondingly low precisions. Moreover, since incomplete factorizations
are by definition inexact, we want to balance the errors introduced by the use of lower precisions
with the error already incurred by the incompleteness of the factorization. This is most straightforward to achieve
for threshold-based factorizations, for which the factorization error is directly controlled by the threshold.
Since threshold-based factorizations do not guarantee any bound on the sparsity of the resulting factors,
more aggressive incomplete factorizations such as level-based or memory-limited ones are also of interest. 
We discuss how these variants can exploit adaptive precision, using a criterion 
based on the magnitude of the largest dropped entry.

The article is structured as follows. In section~\ref{sec:ILUT}, we propose an adaptive precision threshold-based
incomplete factorization and perform an error analysis. In section~\ref{sec:scaling}, we explain how to incorporate
pivoting and scaling within this algorithm. In section~\ref{sec:solves}, we describe possible strategies to perform triangular solves with the resulting
adaptive precision factors.
In section~\ref{sec:IC}, we extend our adaptive precision approach to symmetric positive definite (SPD) matrices and incomplete Cholesky factorizations.
In section~\ref{sec:memlim}, we discuss how this approach can be applied to level-based or memory-limited
factorizations. Throughout these sections, we include representative numerical results to illustrate the discussion. Then, 
in section~\ref{sec:exps}, we perform extensive experiments on a wide range of sparse real-life matrices. Concluding remarks are made in section~\ref{sec:concl}.

Throughout the paper, the unsubscripted norm $\norm{\cdot}$ denotes
the infinity norm $\norm{A}_\infty = \max_i \sum_j \abs{a_{ij}}$.

\subsection{Experimental setting}

In our experiments, we  use sparse matrices taken from the SuiteSparse matrix collection{\footnote{\url{https://sparse.tamu.edu/}}.
In sections~\ref{sec:ILUT-example}, \ref{sec:scaling-example}, \ref{sec:solves-example}, and
\ref{sec:memlim-example}, the  problems listed in Table~\ref{tab:examples} 
are used to illustrate the behavior of the algorithms. 
In section~\ref{sec:exps}, we present results for a much wider test set.

\begin{table}[H]
\caption{Matrices used as illustrative examples.\label{tab:examples}
$\kappa(A)$ is the spectral condition number.
SPD denotes symmetric positive definite. CFD is computational fluid dynamics.}

\begin{tabular}{lrrlcl}
\toprule
Matrix & \multicolumn{1}{c}{$n$} & \multicolumn{1}{c}{Nonzeros} & Sym. & $\kappa(A)$ & Origin \\
\midrule
pores\_2  & 1,224 & 9,613 & Unsym. & $1\times 10^{8\phantom{0}}$ & CFD  \\
radfr1 & 1,048 & 13,299 & Unsym. & $6\times 10^{10}$ & Chemical Process Simulation \\
saylr1 & 238 & 1,128 & Unsym. & $8\times 10^{8\phantom{0}}$ & CFD  \\
nnc666 & 666 & 4,032 & Unsym. & $2\times 10^{10}$ & Nuclear Reactor Modeling \\
bcsstk13 & 2,003 & 83,883 & SPD & $1\times 10^{10}$ & CFD \\
\bottomrule
\end{tabular}
\end{table}

The experiments in sections~\ref{sec:ILUT-example}, \ref{sec:scaling-example}, and \ref{sec:solves-example}
are performed with an ILU factorization implemented in MATLAB. 
The LU factors are used as a left preconditioner for (unrestarted) GMRES.
The experiments in sections~\ref{sec:memlim-example} and \ref{sec:exps} use a modified version of the
Fortran  package {\tt HSL\_MI28} from the HSL mathematical software library\footnote{\url{https://www.hsl.rl.ac.uk/}}. 
For a sparse symmetric positive definite matrix, {\tt HSL\_MI28} computes a memory-limited incomplete Cholesky (IC) factorization
\cite{sctu:2014a,sctu:2014b}. We modify the code to incorporate adaptive precision 
and employ the resulting factors within split-preconditioned CG (PCG).
In our experiments, the stopping tolerance for the relative residual is $10^{-10}$.
We set the maximum number of iterations to $500$ for GMRES and to $2000$ for PCG.

\begin{table}[htbp]
\caption{Floating-point formats used in this work.\label{tab.formats}}
\begin{tabular}{cccccc}
\toprule
Format        & \begin{tabular}[c]{@{}c@{}}Sign\\ bit\end{tabular} &  \begin{tabular}[c]{@{}c@{}}Significand\\ bits\end{tabular} & \begin{tabular}[c]{@{}c@{}}Exponent\\ bits\end{tabular} &\begin{tabular}[c]{@{}c@{}}Range\\ ($f_\mathrm{max}$/$f_\mathrm{min}$)\end{tabular} & Unit Roundoff $\tau$ \\
\midrule
fp64                              & 1    & 52  &  11    & $2^{2046} \approx 10^{616}$  &  $2^{-53} \approx 1 \times 10^{-16}$   \\
rp56                              & 1    & 44  &  11    & $2^{2046} \approx 10^{616}$  & $2^{-45} \approx 3 \times 10^{-14}$   \\
rp48                              & 1    & 36  &  11    & $2^{2046} \approx 10^{616}$  &  $2^{-37} \approx 7 \times 10^{-12}$   \\
rp40                              & 1    & 28  &  11    & $2^{2046} \approx 10^{616}$  &  $2^{-29} \approx 2 \times 10^{-9} $    \\
fp32                              & 1    & 23  &   8    & $2^{254} \approx 10^{76}$      &  $2^{-24} \approx 6 \times 10^{-8} $   \\
rp24                              & 1    & 15  &   8    &  $2^{254} \approx 10^{76}$     &  $2^{-16} \approx 2 \times 10^{-5} $   \\
rp16                     & 1    &  7  &   8    &  $2^{254} \approx 10^{76}$      &  $2^{-8}  \approx  4 \times 10^{-3} $    \\
\bottomrule
\end{tabular}
\end{table}

Our framework and analysis consider a general number $p$ of unspecified precision formats; in practice, we  use the formats
that are
listed in Table~\ref{tab.formats} along with their main characteristics.
The fp64 and fp32 formats are the standard IEEE double and single precisions~\cite{ieee19}. From them we define a set of custom reduced precision
storage formats obtained by reducing the number of bits allocated to the significand: removing one, two, or three bytes from fp64 leads
to the rp56, rp48, and rp40 formats, while removing one or two bytes from fp32 leads to the rp24 and rp16 formats. Note that
rp16 is equivalent to the bfloat16 format~\cite{inte18}, although in this work we only use it as a storage format
and do not require hardware support for bfloat16 arithmetic.
In total we have seven different precision formats, with fine variations on the unit roundoff
that allow us to flexibly control the accuracy of the factor entries.
These reduced precision formats have been successfully used in other contexts
\cite{agjm, ajlm26, bmv, gjmm24a}. 
Note that the formats used in this work have either 8 or 11 exponent bits, and thus provide
a sufficiently wide range to prevent underflows and overflows in practice.
Exploring the use of reduced exponent formats is an interesting perspective for future work.

\section{Adaptive precision threshold-based factorization} \label{sec:ILUT}

\subsection{Algorithm}
Algorithm~\ref{alg.ILUT} outlines the proposed approach: an adaptive precision incomplete LU (hereafter \APILU) factorization.
We denote the target accuracy by $\epsabs$: as we explain below, it is used to decide
which factor entries are dropped, and the precision used
to store the entries that are kept.
For simplicity of notation, the algorithm
is outlined for dense matrices; it can be modified for sparse matrices by looping the indices over  nonzero entries of the factors.

\begin{algorithm}[htbp]
\caption{Adaptive precision ILU (\APILU)\label{alg.ILUT}}
\textbf{Input:}  \begin{tabular}[t]{l}$A\in\R^{n\times n}$, target accuracy $\epsabs$, storage precisions $\{\tau_1,\ldots,\tau_p\}$ with \\$\tau_1 < \ldots < \tau_p$, compute precision $\zeta$  \end{tabular} \\
\textbf{Output:} $\Lh\Uh\approx A$

\begin{algorithmic}[1]
\FOR{$k=1\colon n$}
\FOR{$j=k\colon n$}
\STATE $u_{kj} = a_{kj} - \sum_{i=1}^{k-1} \lh_{ki}\uh_{ij}$ in precision $\zeta$ \label{line.sum1}
\STATE $\tau = \max\left\{\tau \in \{\tau_1,\ldots,\tau_p\} \colon \abs{u_{kj}} \le \epsabs/\tau \right\}$ \label{line.tauU}
\STATE $\uh_{kj} = \fl_\tau(u_{kj})$ \label{line.fl1}
\ENDFOR
\FOR{$i=k+1\colon n$}
\STATE $\ell_{ik} = \big(a_{ik} - \sum_{j=1}^{k-1} \lh_{ij}\uh_{jk}\big)/\uh_{kk}$ in precision $\zeta$ \label{line.sum2}
\STATE $\tau = \max\left\{\tau \in \{\tau_1,\ldots,\tau_p\} \colon \abs{\ell_{ik}\uh_{kk}} \le \epsabs/\tau \right\}$ \label{line.tauL}
\STATE $\lh_{ik} = \fl_\tau(\ell_{ik})$ \label{line.fl2}
\ENDFOR
\ENDFOR
\end{algorithmic}
\end{algorithm}

Algorithm~\ref{alg.ILUT} is an adaptive precision method: it uses
$p$ different floating-point precision formats with unit roundoffs $\{\tau_1,\ldots,\tau_p\}$ satisfying $\tau_1 < \ldots < \tau_p$.
We define $\fl_\tau(\cdot)$ to be the operator that rounds to precision $\tau$: for any $\alpha\in\R$ such that 
$\abs{\alpha}$ belongs to the representable range $[f_\mathrm{min}, f_\mathrm{max}]$ of the precision format,
we have~\cite[Theorem~2.2]{high:ASNA2}:
\begin{equation}\label{eq.std}
\fl_\tau(\alpha) = \alpha(1+\delta), \quad \abs{\delta} \le \tau.
\end{equation}
Note that the usual dropping strategy (replacing a number by zero) can be incorporated in this framework by taking the last ``precision'' to be $\tau_p=1$, 
and defining $\fl_1(\cdot) = 0$ (which satisfies \eqref{eq.std} with $\delta=-1$).

The precisions $\{\tau_1,\ldots,\tau_p\}$ are used for \textit{storage only}. We  employ a \textit{compute precision} 
with unit roundoff $\zeta$ for performing the computations on lines \ref{line.sum1} and \ref{line.sum2}.
We use the standard model of floating-point arithmetic~\cite[(2.4)]{high:ASNA2}
\begin{equation}\label{eq.std-2}
\fl_\zeta(\alpha \op \beta) = (\alpha\op\beta)(1+\delta), \quad \abs{\delta} \le \zeta, \quad \op\in\{+, -, \times, \div\}.
\end{equation}
In practice, we use fp64 arithmetic for the computations; 
thus, the typical setting of interest 
is $\zeta \ll \epsabs$ and we can expect
the accuracy of the LU factors to be mainly determined by the
choice of storage precisions.
If the conditions on lines~\ref{line.tauU} and \ref{line.tauL} are not satisfied by any of
the storage precisions, the entry is retained in the compute precision, so that no rounding
error occurs.

\subsection{Analysis} \label{sec:ILUT-analysis}

Algorithm~\ref{alg.ILUT} is based on the key observation 
that we can store each factor entry in a precision with a unit roundoff $\tau$
inversely proportional to the magnitude of the entry.
The following theorem proves that Algorithm~\ref{alg.ILUT} is backward stable
with respect to the target accuracy $\epsabs$.

\begin{theorem}\label{thm.ILUT}
If Algorithm~\ref{alg.ILUT} applied to $A\in\R^{n\times n}$ runs to completion, the computed LU factors satisfy
\begin{equation} \label{eq.bnd}
A+E = \Lh\Uh,\quad \abs{E} \le (1 + \gamma_{\ceps})\epsabs J + \gamma_{\ceps} \abs{\Lh}\abs{\Uh},
\end{equation}
where $J$ is the matrix of ones,
$\ceps \le n$ is the maximum number of nonzero entries in any row or column of $\Lh$ and $\Uh$,
and $\gamma_{\ceps} = \ceps\,\zeta / (1 - \ceps \,\zeta)$.
\end{theorem}

\begin{proof}
  We  use the standard notations~\cite{high:ASNA2} $\theta_{c}$ and $\gamma_{c}$ for analyzing cumulative
  rounding errors, where given $\abs{\delta_{i}} \le \zeta$, $\rho_{i} = \pm 1$, and assuming $c\,\zeta < 1$, we define
  \[
    (1 + \theta_{c}) = \prod_{i=1}^{c}(1 + \delta_{i})^{\rho_{i}},
    \quad \abs{\theta_{c}} \leq \frac{c \,\zeta}{1 - c \,\zeta} = \gamma_{c}.
  \]
  For the upper triangle, $a_{kj}, 1 \leq k \leq j \leq n$, we have, from line \ref{line.sum1},
  \[
  u_{kj} = a_{kj} - \sum_{i=1}^{k-1} \lh_{ki} \uh_{ij}.
  \]
  Using the standard analysis of inner products in LU factorization~\cite[Chap.~9]{high:ASNA2},
  the computed $\ut_{kj}$ in precision $\zeta$ satisfies
  \begin{equation} \label{eq.thm1.vhat}
    \ut_{kj}(1+\theta_{\ceps}^{(k)}) = a_{kj}
    - \sum_{i=1}^{k-1} \lh_{ki} \uh_{ij} (1+\theta_{\ceps}^{(i)}).
  \end{equation}
  Then, line~\ref{line.fl1} stores the computed $\ut_{kj}$ in its selected storage precision $\tau$, yielding
  \begin{equation} \label{eq.thm1.uhat}
    \uh_{kj} = \fl_{\tau}(\ut_{kj}) = (1 + \delta)\ut_{kj} = \ut_{kj} + e_{kj},
  \end{equation}
  where $\abs{\delta} \le \tau$ and $\tau$ is chosen such that
  $\abs{e_{kj}} = \abs{\delta \ut_{kj}} \leq \tau \abs{\ut_{kj}} \le \epsabs$.
  Substituting into \eqref{eq.thm1.vhat}, we obtain
\[
    (\uh_{kj}-e_{kj})(1+\theta_{\ceps}^{(k)}) = a_{kj}
    - \sum_{i=1}^{k-1} \lh_{ki} \uh_{ij} (1+\theta_{\ceps}^{(i)}),
\]
and since $\lh_{kk} = 1$, this can be rewritten to get a backward error result
  \begin{equation}\label{eq.bndU}
  \begin{aligned}
  \sum_{i=1}^{k}\lh_{ki} \uh_{ij} &= a_{kj} + (1 + \theta_{\ceps}^{(k)})e_{kj} - \sum_{i=1}^{k}\lh_{ki} \uh_{ij} \theta_{\ceps}^{(i)} \\
  &= a_{kj} + \delta a_{kj}, \quad \abs{\delta a_{kj}} \le (1 + \gamma_{\ceps})\epsabs + \gamma_{\ceps} \sum_{i=1}^{k} \abs{\lh_{ki}} \abs{\uh_{ij}}. 
  \end{aligned}
  \end{equation}

  For the lower triangle, $a_{ik}, 1 \leq k < i \leq n$, we can similarly prove that
  \[
  \sum_{j=1}^{k}\lh_{ij} \uh_{jk} = a_{ik} + \delta a_{ik}, \quad \abs{\delta a_{ik}} \le (1 + \gamma_{\ceps})\epsabs + \gamma_{\ceps} \sum_{j=1}^{k} \abs{\lh_{ij}} \abs{\uh_{jk}}.
  \]

  The proof is completed by writing the errors in matrix form.
\end{proof}

\begin{remark}\label{rem.epschoice}
Note that the first term, $(1 + \gamma_{\ceps})\epsabs J$, in bound \eqref{eq.bnd} is an absolute quantity. 
A relative bound can
be obtained by scaling $\epsabs$ by an appropriate quantity.
For example, a relative normwise backward error bound of order $\epsrel$ can be obtained by taking $\epsabs =\epsrel\norm{A}$.
Alternatively, if we want a componentwise bound, $\epsabs$ needs to be made dependent on the indices, such as $\varepsilon_{ik}=\epsrel\abs{a_{ik}}$
or perhaps $\varepsilon_{ik} = \epsrel\sum_{j=1}^{k} \abs{\lh_{ij}\uh_{jk}}$, the latter choice leading to a bound $\abs{E} \lesssim \epsrel\abs{\Lh}\abs{\Uh}$
that resembles the usual backward error bound in floating-point arithmetic with unit roundoff $\epsrel$.
Throughout the rest of this paper, we focus on normwise bounds and fix $\epsabs=\epsrel\norm{A}$ for all entries.
\end{remark}

\begin{remark}\label{rem.facto-time}
With a left-looking implementation of Algorithm~\ref{alg.ILUT}, the computations on lines~\ref{line.sum1} and \ref{line.sum2}
require accessing the previously computed entries of the factors, which are already stored in adaptive precision. Therefore, in a regime
where the computation is memory-bound, in addition to reducing storage, Algorithm~\ref{alg.ILUT} may improve the computational time,
although a precise assessment of this requires a dedicated performance study and is outside the scope of this numerical study.
\end{remark}

\begin{remark}\label{rem.AP-compute}
It is natural to ask whether the computations on lines~\ref{line.sum1} and~\ref{line.sum2} can be performed in adaptive precision,
that is, whether it is meaningful to replace the uniform compute precision $\zeta$ with multiple precisions $\zeta_1, \ldots, \zeta_q$. 
In principle, this is possible: Theorem~\ref{thm.ILUT} shows that the error incurred on an entry $\uh_{kj}$ ($k \le j$)
is proportional to $\sum_{i=1}^k \abs{\lh_{ki}}\abs{\uh_{ij}}$, so this entry could be computed in precision 
$\zeta = \max\{\zeta\in\{\zeta_1,\ldots,\zeta_q\} : \sum_{i=1}^k \abs{\lh_{ki}}\abs{\uh_{ij}} \le \epsabs/\zeta\}$.
Two practical difficulties arise.
First, $\sum_{i=1}^k \abs{\lh_{ki}}\abs{\uh_{ij}}$ can be substantially larger than $\abs{\uh_{kj}}$, meaning that 
$\uh_{kj}$ may need to be computed in a precision with unit roundoff $\zeta \ll \tau$ much smaller
than that of its storage precision.
Second,  $\sum_{i=1}^k \abs{\lh_{ki}}\abs{\uh_{ij}}$ is unknown and computing it requires about as many operations
as computing $\uh_{kj}$. One possible strategy is to compute $\sum_{i=1}^k \abs{\lh_{ki}}\abs{\uh_{ij}}$
in the lowest precision $\zeta_q$ to estimate its order of magnitude, and then
compute $\uh_{kj}$ in the required precision. This type of look-ahead
approach has been successfully employed in other contexts~\cite{efmr,bgyt26}.
Alternatively, a max-plus approach~\cite{hsth:2018,hoti17} can potentially be used to estimate the magnitude of the entries of the LU factors.
These ideas are interesting but difficult to put into practice; in this paper, we focus on adaptive precision storage
and leave the computations in uniform, high precision.
\end{remark}

\subsection{Numerical example}
\label{sec:ILUT-example}
\begin{figure}[ht]
    \centering
    \begin{subfigure}[c]{0.49\textwidth}
        \centering
        \includegraphics[width=\linewidth]{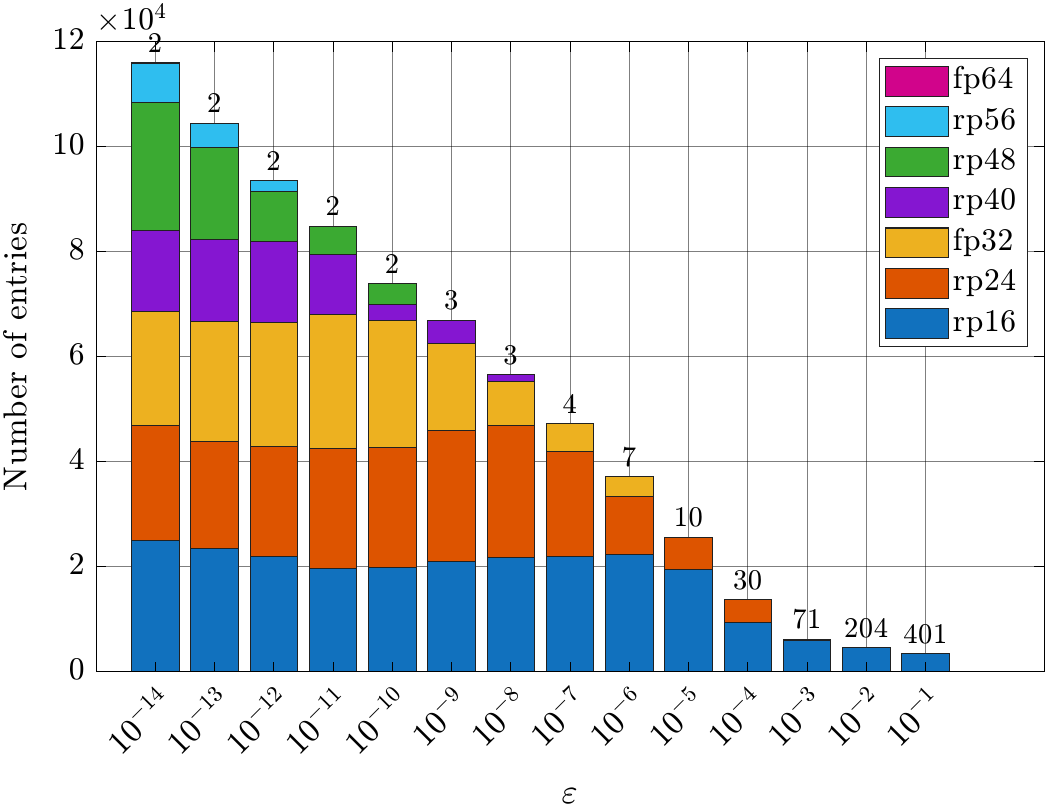}
        \caption{Distribution of entries across precisions}
        \label{fig:pores_2_distrib}
    \end{subfigure}
    \hfill
    \begin{subfigure}[c]{0.49\textwidth}
        \centering
        \includegraphics[width=\linewidth]{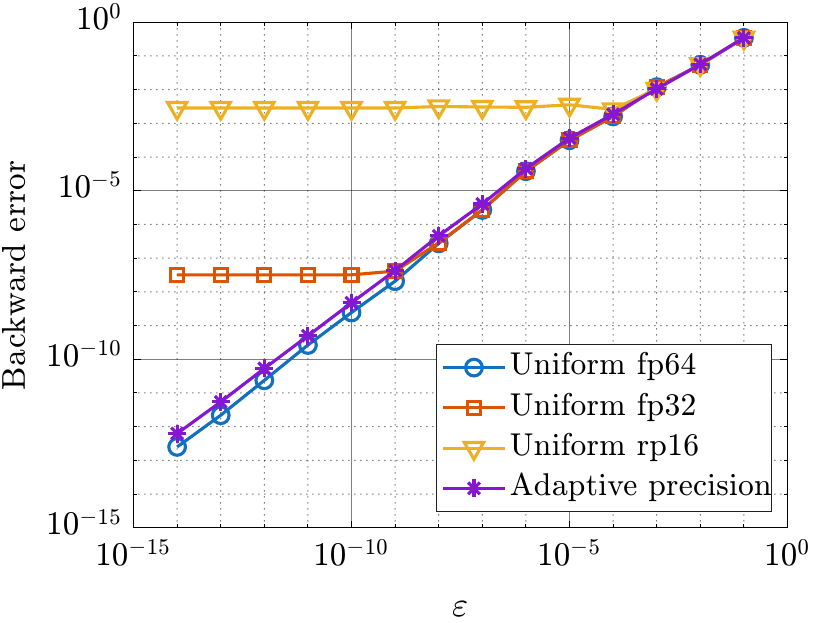}
        \caption{Error vs $\epsrel$}
        \label{fig:pores_2_berr_eps}
    \end{subfigure}
    \begin{subfigure}[c]{0.49\textwidth}
        \centering
        \includegraphics[width=\linewidth]{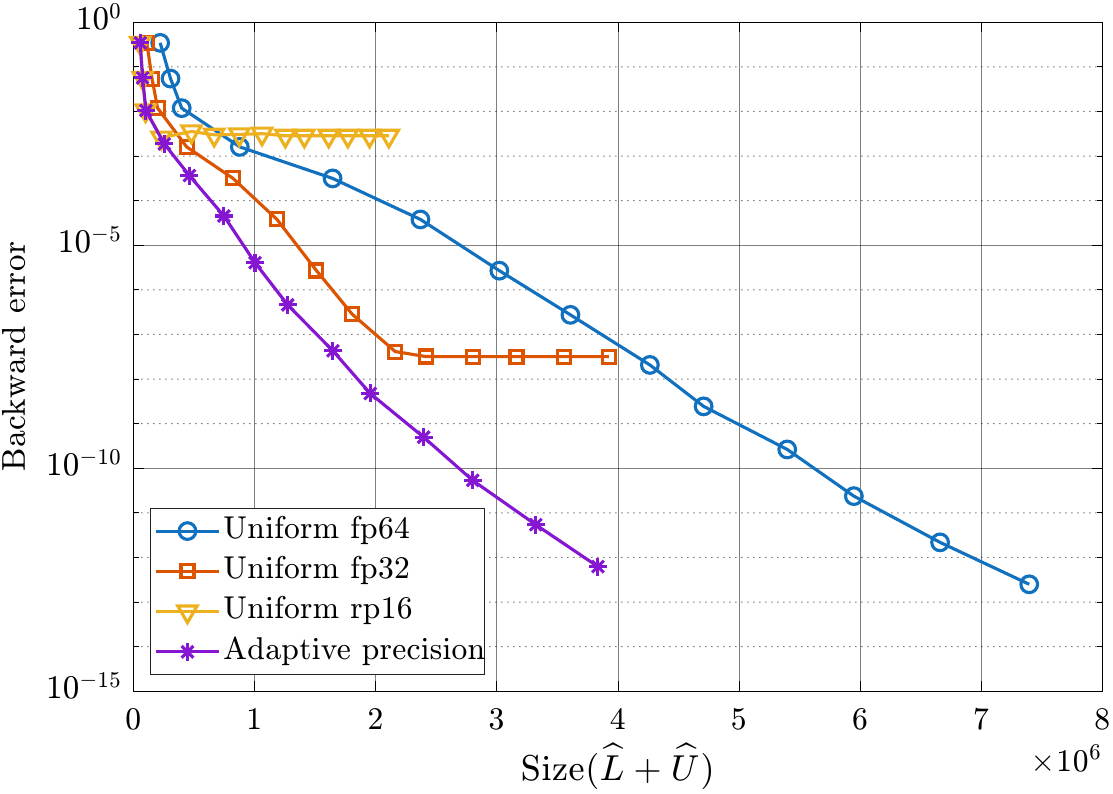}
        \caption{Error vs size$(\Lh+\Uh)$}
        \label{fig:pores_2_berr_byte}
    \end{subfigure}
    \hfill
    \begin{subfigure}[c]{0.49\textwidth}
        \centering
        \includegraphics[width=\linewidth]{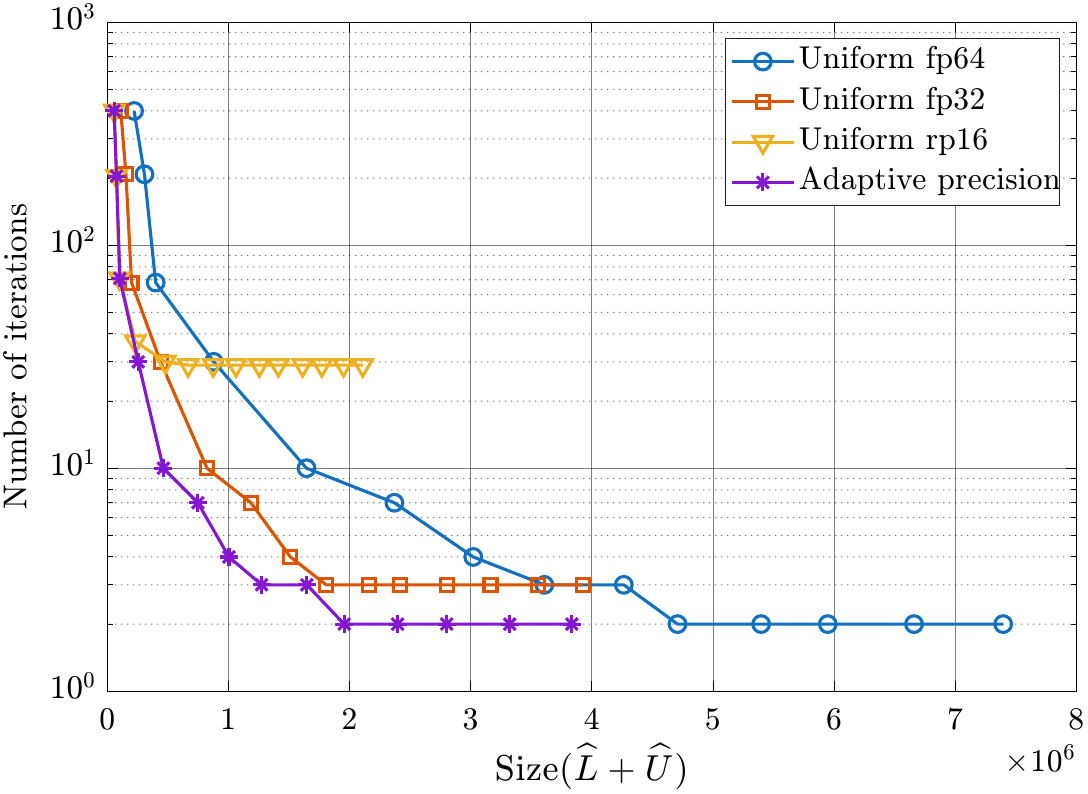}
        \caption{GMRES iterations vs size$(\Lh+\Uh)$}
        \label{fig:pores_2_iter_byte}
    \end{subfigure}
    \caption{Results for Algorithm~\ref{alg.ILUT} on matrix pores\_2.}
    \label{fig:pores_2}
\end{figure}

The behavior of Algorithm~\ref{alg.ILUT} 
depends on the distribution of the magnitudes
of the entries of the factors, and thus on $A$. The goal of this section
is to shed some light on this using a numerical example;
a more quantitative assessment on a wide range of matrices is given in section~\ref{sec:exps}.

Figure~\ref{fig:pores_2} shows the result of running Algorithm~\ref{alg.ILUT} (\APILU)
on the unsymmetric matrix pores\_2.
Figure~\ref{fig:pores_2_distrib} plots the distribution of the LU entries across
the seven precision formats, depending on $\epsrel$; the height of the bars 
corresponds to the number of LU entries and hence the amount of dropping. The numbers on top
of the bars indicate the number of GMRES iterations plotted in Figure~\ref{fig:pores_2_iter_byte}
(see description below).
This first plot illustrates that, as $\epsrel$ increases, we drop more entries
(as in uniform precision ILU) while also switching more and more of the retained entries to lower and lower precisions.
For example, at $\epsrel=10^{-14}$, despite the factorization being nearly exact, we store
the majority of the entries in rp48 or lower precisions; fewer than $0.1\%$ of the entries need to be stored in the highest fp64 precision. 
At $\epsrel=10^{-5}$, most entries can be stored in rp16.

Figure~\ref{fig:pores_2_berr_eps} plots the relative normwise backward error $\norm{A-\Lh\Uh}/\norm{A}$
for \APILU as well as its uniform precision baselines, for which
all entries are stored in the same uniform precision (fp64, fp32, or rp16).
We see that for the uniform fp64 ILU, the error is proportional
to $\epsrel$. As expected from Theorem~\ref{thm.ILUT}, this is also true for \APILU,
which achieves an error that is comparable to uniform fp64 for small values of $\epsrel$ and
indistinguishable from it as $\epsrel$ grows. Hence the reduction of storage achieved by adaptive precision is essentially
obtained without any significant loss of accuracy; this illustrates that \APILU
successfully balances the introduced rounding errors with the dropping errors, so that the effect
of the former is hidden by that of the latter.
In contrast, lower precision uniform baselines, fp32 and rp16, cannot achieve an error below
the unit roundoff of their precision; their backward errors are the same as for \APILU only when
$\epsrel \ge 10^{-9}$ and $\epsrel \ge 10^{-4}$, respectively.

Figure~\ref{fig:pores_2_berr_byte} plots the same backward error, but against the storage
required by the LU factors on the $x$-axis; the storage is measured by accounting for the 
number of bits required by each factor entry depending on its precision.
Here and for our other MATLAB experiments, we report the storage cost
of the floating-point values and ignore the cost of storing the integer indices required for sparse implementations;
the storage cost of these indices will be evaluated in section~\ref{sec:memlim} with the \texttt{HSL\_MI28} experiments.
This figure shows that \APILU can achieve significant storage reductions 
compared with uniform precision baselines for a given prescribed error; equivalently, for a fixed
memory budget, \APILU can significantly improve accuracy. For example, to achieve an error below
$10^{-10}$, uniform fp64 ILU requires about 6~Mb, whereas \APILU requires less than 3~Mb.
These reductions become even more significant if larger errors are allowed.
For example, if an error of $10^{-5}$ is acceptable, the storage for uniform fp64 ILU reduces by a factor of $2\times$ (from 6 to 3~Mb), 
whereas that for \APILU reduces by a factor $3\times$ (from 3 to 1~Mb).
The comparison with lower precision uniform baselines leads to the same conclusions; they can
reduce the storage but are limited to errors above a certain threshold. \APILU consistently achieves
the best storage--accuracy tradeoff.

Figure~\ref{fig:pores_2_iter_byte} plots the number of iterations
required for GMRES to converge when the \APILU factors are used as the preconditioner.
Comparing \APILU with the fp64 and fp32 uniform precision baselines shows that \APILU
achieves a significantly better storage--iteration tradeoff: \APILU can converge in fewer than 10 iterations while
reducing the storage of 1.6 and 0.8 Mb needed by uniform fp64 and fp32, respectively,
to 0.5 Mb. 
Uniform rp16 can achieve low storage, but always requires at least 30 iterations to converge.
It is interesting to compare \APILU and uniform rp16 at $\epsrel=10^{-5}$. As shown in Figure~\ref{fig:pores_2_distrib}, 
\APILU stores most of the entries in rp16, with about 20\% stored in rp24.
This has little effect on the total storage but reduces the
iteration count from 30 to 10.

\begin{figure}
    \centering
    \begin{subfigure}[t]{0.49\textwidth}
        \vspace{0pt}
        \centering
        \includegraphics[width=\linewidth]{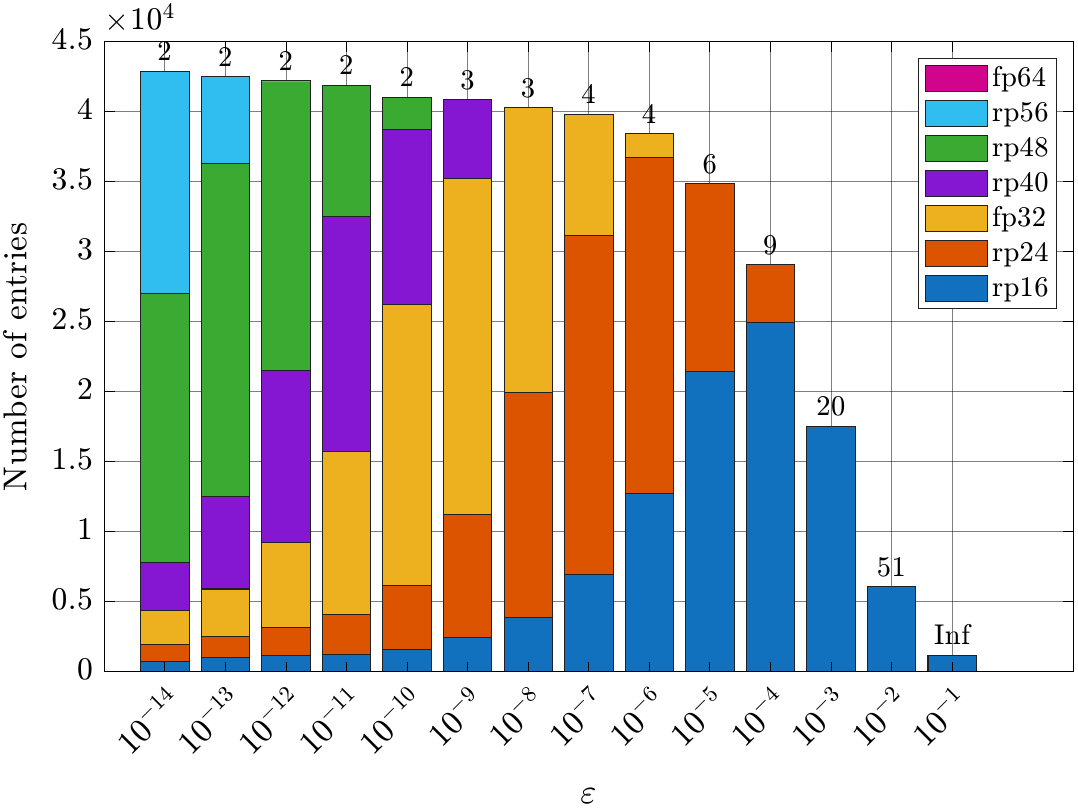}
        \caption{Distribution of entries across precisions. ``Inf'' indicates that GMRES failed to converge.}
        \label{fig:radfr1_distrib}
    \end{subfigure}
    \hfill
    \begin{subfigure}[t]{0.49\textwidth}
        \vspace{0pt}
        \centering
        \includegraphics[width=\linewidth]{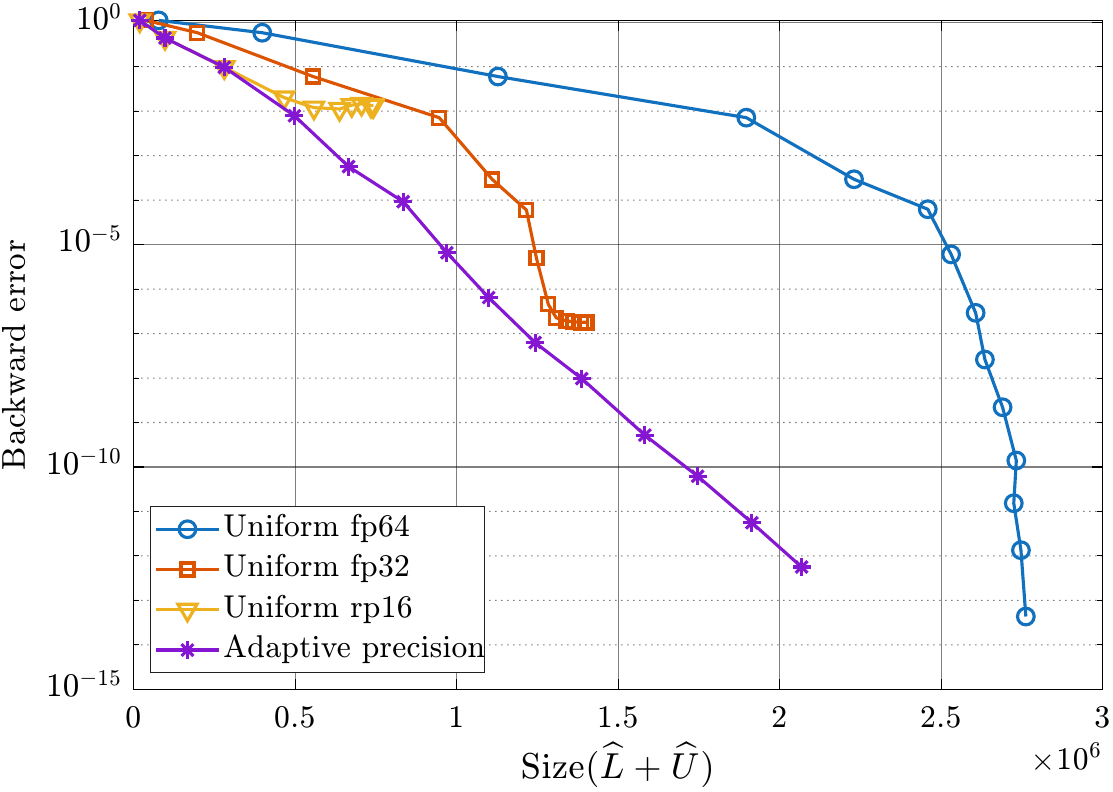}
        \caption{Error vs size$(\Lh+\Uh)$}
        \label{fig:radfr1_berr_byte}
    \end{subfigure}
    \caption{Results for Algorithm~\ref{alg.ILUT} on matrix radfr1.}
    \label{fig:radfr1}
\end{figure}

We conclude this section by reporting results for matrix radfr1 in Figure~\ref{fig:radfr1}.
This illustrates an interesting feature of \APILU: it can perform well
even on matrices where classical (uniform precision) ILU does not. For this matrix,
few entries are dropped for $\epsrel \le 10^{-6}$, which leads to very small storage
reductions  with uniform fp64 ILU.
In contrast, \APILU
compensates for the relatively small amount of dropping allowed by this matrix with adaptive precision storage. As shown
in Figure~\ref{fig:radfr1_distrib}, even though about the same number of entries are retained as $\epsrel$ increases,
their precision distribution significantly shifts towards the lower precisions. At $\epsrel=10^{-14}$, most of the entries
must be stored in rp56 and rp48, leading to a storage cost of 2.1~Mb, whereas at $\epsrel=10^{-6}$, most of the entries
can be stored in rp24 and rp16, leading to a reduced storage cost of 0.8~Mb.
Despite this reduction, the  preconditioner is of high quality, with only 4 iterations
required for convergence.

\section{Incorporating scaling and pivoting}
\label{sec:scaling}

ILU factorizations  usually incorporate scaling and pivoting to improve the quality
of the preconditioner. It is therefore an important practical concern to extend Algorithm~\ref{alg.ILUT}
and its analysis to incorporate these techniques.

\subsection{Pivoting}
Incomplete factorizations can break down because of zero (or near-zero) pivots (this can occur regardless of the precision used---even in exact arithmetic). 
Pivoting strategies, adapted from sparse direct solvers, address this by reordering
rows and/or columns so that an entry of large magnitude is permuted onto the $j$th diagonal position
before  elimination step $j$ begins. 
The resulting preconditioners have been systematically analyzed and popularized by
Saad~\cite{Saad1994,saad:2003}; other key references include \cite{BenziSzyldVanDuin1999,bosa:2006,ChowSaad1997,lish:2011}.

It is straightforward to account for pivoting in our analysis: computing a factorization  $PAQ \approx \Lh\Uh$ is  equivalent (both in exact and finite precision
arithmetic) to computing an (unpivoted) ILU 
factorization of the prepermuted matrix $PAQ$. Therefore, the analysis and the resulting criteria to adapt the precisions are unchanged.

\subsection{Scaling}

We next discuss scaling. Let $B= \Lambda^{-1} A \Theta^{-1}$ be the scaled matrix, where $\Lambda=\diag(\lambda_i)$ and $\Theta=\diag(\theta_i)$ are diagonal matrices.
We distinguish two possible strategies.
A first approach is to directly
apply Algorithm~\ref{alg.ILUT} to $B$; uniform precision incomplete factorizations use this approach for choosing which entries to drop. 
By Theorem~\ref{thm.ILUT}, and neglecting terms of order $\zeta$, this yields
\begin{align}
& B + E = \Lh\Uh,\quad \abs{E} \le \epsabs J \\
\Leftrightarrow\quad & A + \Lambda E\Theta = \Lambda\Lh\Uh\Theta,\quad \abs{\Lambda E\Theta} \le \epsabs \abs{\Lambda} J\abs{\Theta}, \label{eq.LJT}
\end{align}
where $J$ is the matrix of ones.
The term $\abs{\Lambda}J\abs{\Theta}$ shows that we obtain an error bound that depends on the scaling.
In the common case where scaling is used to equilibrate the rows and/or columns of the matrix, 
\eqref{eq.LJT} essentially reduces to a relative rowwise or columnwise bound.
For example, if $\lambda_i = 1$ and $\theta_j = \norm{a_{:j}}$, where $a_{:j}$ denotes the $j$th column of $A$, 
then \eqref{eq.LJT} yields a componentwise error $\varepsilon_{ij} = \epsabs\norm{a_{:j}}$.

However, if only a normwise bound is needed, \eqref{eq.LJT}
is unsatisfactory. We would like to have instead a uniform error bound of the form $\abs{\Lambda E\Theta} \le \epsabs J$. 
To achieve this, we must take into account
the scaling matrices within the precision choice criterion in Algorithm~\ref{alg.ILUT}.
Specifically, the equation analogous to \eqref{eq.bndU}
when applied to $B$ instead of $A$ (neglecting again terms of order $\zeta$) is
\begin{align}
& b_{kj} + e_{kj} = \sum_{i=1}^{k} \lh_{ki}\uh_{ij}, \quad \abs{e_{kj}}\le\tau\abs{u_{kj}} \\
\Leftrightarrow\quad 
& a_{kj} + \lambda_{k}e_{kj}\theta_{j} = \sum_{i=1}^{k} \lambda_{k}\lh_{ki}\uh_{ij}\theta_{j}, \quad \abs{\lambda_{k} e_{kj} \theta_{j}}\le\tau\abs{\lambda_{k}u_{kj}\theta_{j}}.
\label{eq.bnd-scaled}
\end{align}
Thus line~\ref{line.tauU} in Algorithm~\ref{alg.ILUT} should be modified to
\begin{equation}\label{eq.unscaled-modif1}
\tau = \max\left\{\tau \in \{\tau_1,\ldots,\tau_p\} \colon \abs{\lambda_{k} u_{kj} \theta_{j}} \le \epsabs/\tau \right\}
\end{equation}
By the same reasoning, line~\ref{line.tauL} should be modified to
\begin{equation}\label{eq.unscaled-modif2}
\tau = \max\left\{\tau \in \{\tau_1,\ldots,\tau_p\} \colon \abs{\lambda_{i}\ell_{ik}\uh_{kk}\theta_{k}} \le \epsabs/\tau \right\}
\end{equation}
With these modifications \eqref{eq.bnd-scaled} reduces to a uniform bound $\epsabs J$.

\subsection{Numerical example}
\label{sec:scaling-example}
Figure~\ref{fig:scaling} compares the two approaches described in the previous section using matrix saylr1.
The ``scaled approach'' applies \APILU to $B=\Lambda^{-1}A\Theta^{-1}$ with $\epsabs=\epsrel\norm{B}$, whereas the ``unscaled approach'' also works
on $B$, but modifies the criterion to choose the precisions (including dropping) as in \eqref{eq.unscaled-modif1}--\eqref{eq.unscaled-modif2}, with $\epsabs=\epsrel\norm{A}$,
that is, it \textit{undoes the scaling} for the dropping and precision choice.
For this example, we 
use $\lambda_i=\max_j \abs{a_{ij}}$ and $\theta_j = \max_i \abs{a_{ij}}$.

\begin{figure}
    \centering
    \begin{subfigure}{0.49\textwidth}
        \centering
        \includegraphics[width=\linewidth]{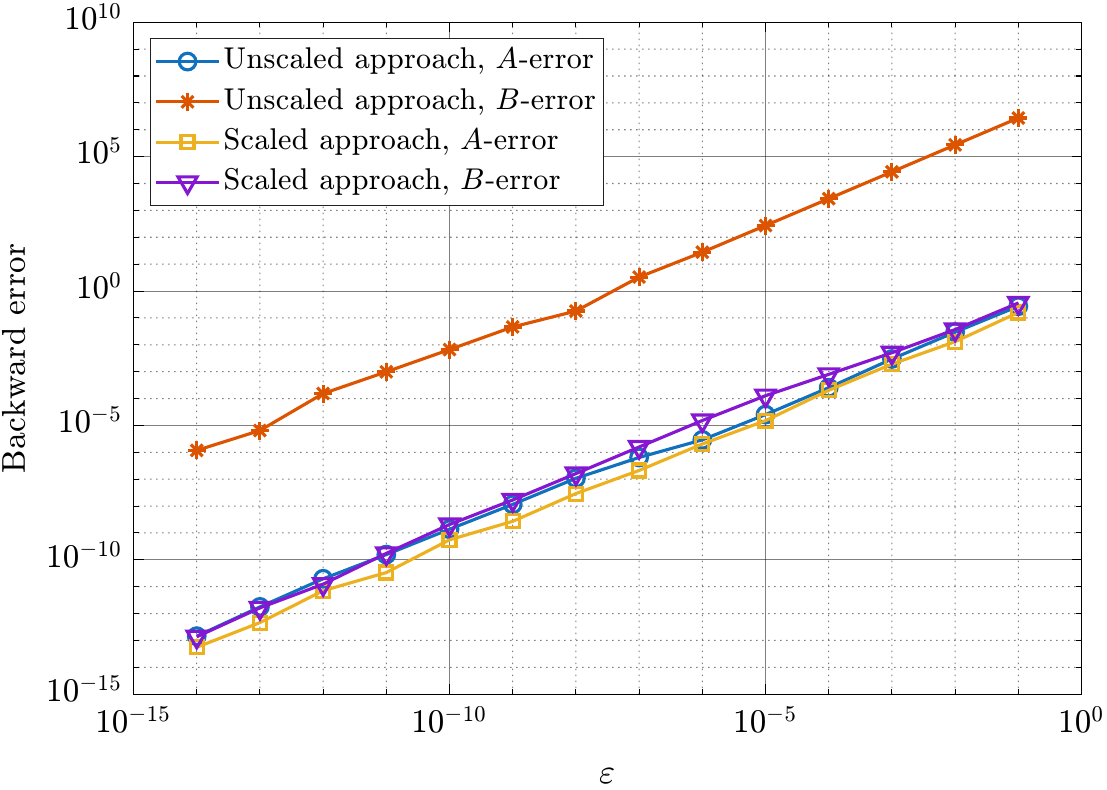}
        \caption{Error vs $\epsrel$}
        \label{fig:scaling_saylr1_berr_eps}
    \end{subfigure}
    \hfill
    \begin{subfigure}{0.49\textwidth}
        \centering
        \includegraphics[width=\linewidth]{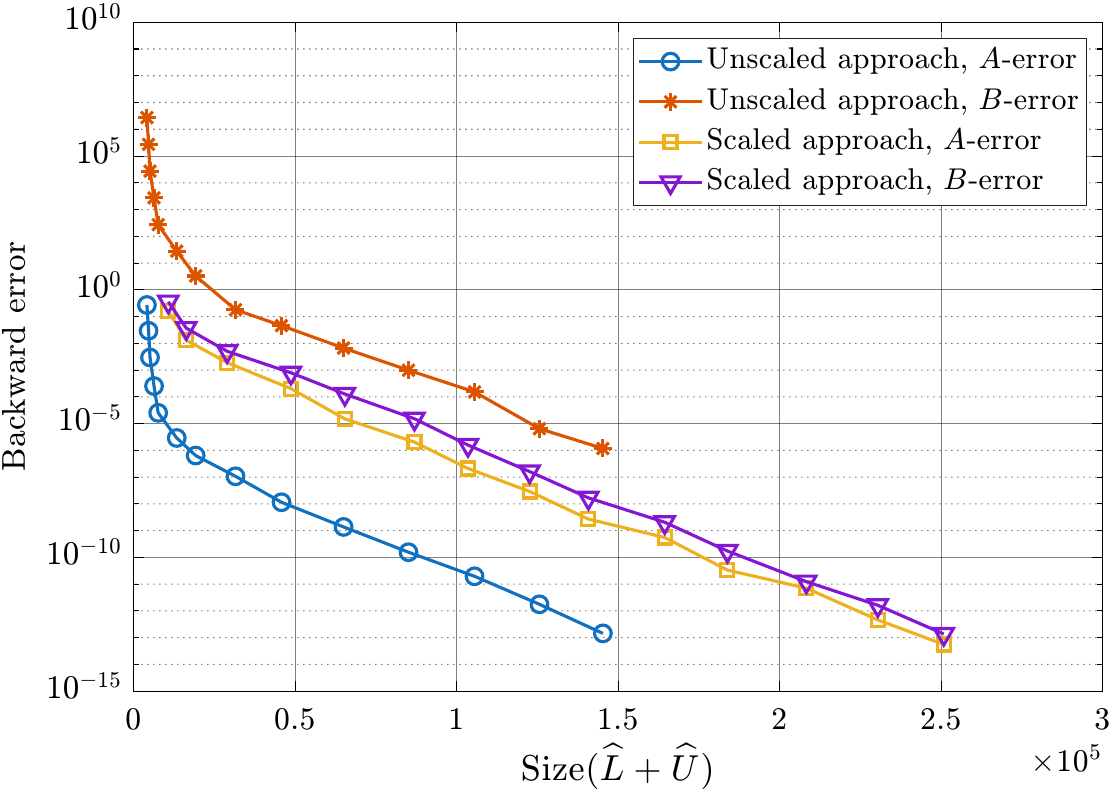}
        \caption{Error vs size$(\Lh+\Uh)$}
        \label{fig:scaling_saylr1_berr_byte}
    \end{subfigure}
    \begin{subfigure}{0.49\textwidth}
        \centering
        \includegraphics[width=\linewidth]{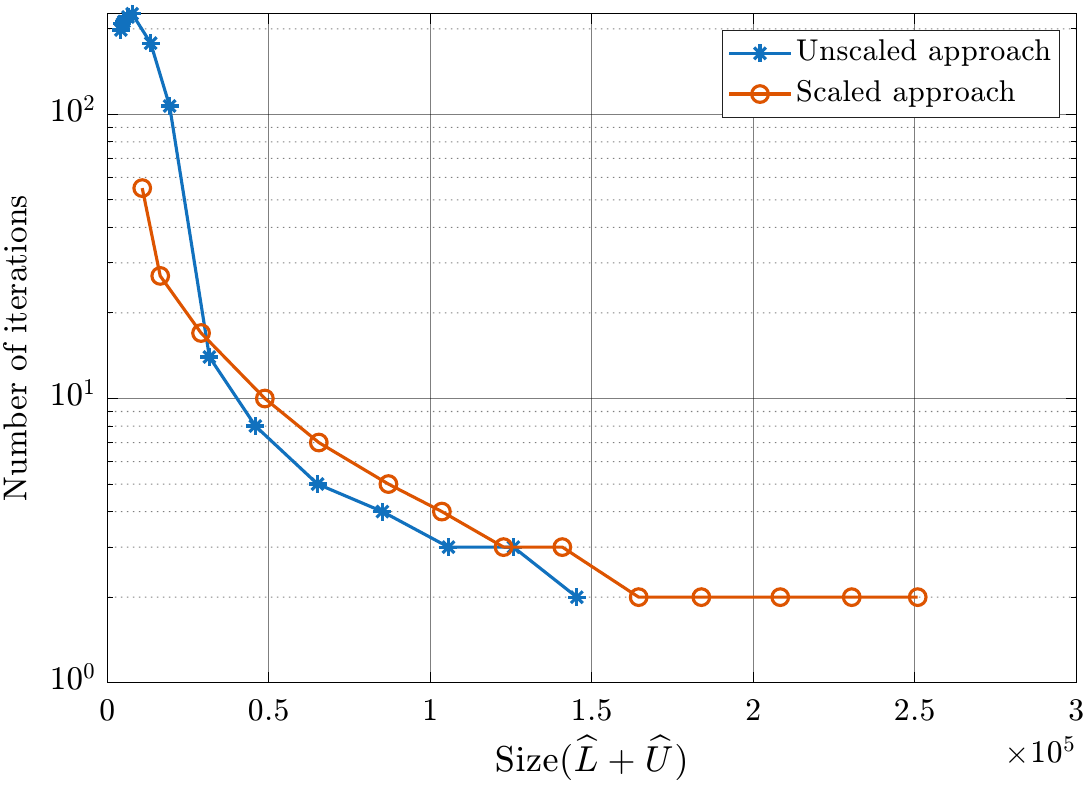}
        \caption{GMRES iterations vs size$(\Lh+\Uh)$}
        \label{fig:scaling_saylr1_iter_byte}
    \end{subfigure}
    \caption{Comparison between scaled and unscaled approaches on matrix saylr1.}
    \label{fig:scaling}
\end{figure}

Figure~\ref{fig:scaling_saylr1_berr_eps} plots two variants of the relative normwise backward error:
the ``$B$-error'', which is the error on $B$:
\begin{equation}
B\text{-error} = \frac{\norm{B-\Lh\Uh}}{\norm{B}} = \frac{\norm{\Lambda^{-1}A\Theta^{-1} - \Lh\Uh}}{\norm{\Lambda^{-1}A\Theta^{-1}}},
\end{equation}
and the ``$A$-error'', which is the error on $A$:
\begin{equation}
A\text{-error} = \frac{\norm{A-\Lambda\Lh\Uh\Theta}}{\norm{A}}.
\end{equation}
Figure~\ref{fig:scaling_saylr1_berr_byte} plots the same errors against the factor storage cost on the $x$-axis.
For the scaled approach, the $A$- and $B$-errors
are similar, but are quite different for the unscaled approach. The $A$-error for the unscaled
approach is similar to both errors for the scaled approach with the same $\epsrel$ (see Figure~\ref{fig:scaling_saylr1_berr_eps}). 
However, the storage cost for the unscaled approach is significantly lower, which means that when measuring accuracy in
the $A$-error, the unscaled approach achieves a better storage--accuracy tradeoff than the scaled approach (see Figure~\ref{fig:scaling_saylr1_berr_byte}).
This is expected: both approaches provide an $A$-error
of order $\epsrel$, but the unscaled approach can have more opportunities to use lower precisions (and dropping)
because the criteria 
\eqref{eq.unscaled-modif1}--\eqref{eq.unscaled-modif2} consider the unscaled entries, which should
exhibit wider variations in magnitude.
However, the unscaled approach leads to a $B$-error
that is much larger; this is because it introduces large perturbations
relative to the scaled entries of $B$. The $B$-error can be seen as a form
of rowwise and columnwise error, whereas the $A$-error is only a globally normwise error.

Figure~\ref{fig:scaling_saylr1_iter_byte} shows that there is no universal
answer as to which error metric is the most relevant when using \APILU as a preconditioner. It plots the number of iterations required for GMRES to converge against the storage cost.
For small $\epsrel$, the unscaled approach achieves a better storage--iteration tradeoff than the scaled one,
 whereas the converse is true for large $\epsrel$. 
 We have performed this same comparison on a wide range of matrices, and which of the scaled
 or unscaled approaches achieves the better tradeoff is heavily matrix dependent.
Some matrices may benefit from a small scaled (rowwise/columnwise) error to converge rapidly, whereas
for other matrices,  a small globally normwise error is sufficient.
In general, we found the scaled approach to be more robust and consistent overall. 
The unscaled approach can give greater storage savings,  but its aggressive dropping strategy sometimes prevents convergence unless
$\epsrel$ is very small, and its behavior is thus more sensitive to the choice of $\epsrel$. 
Therefore, throughout the rest of this study, we focus on the scaled approach
(note that the previously presented examples in section~\ref{sec:ILUT-example} used the scaled approach).

\section{Adaptive precision triangular solves}
\label{sec:solves}

Once the incomplete factorization has been computed, each application of
the resulting preconditioner involves performing triangular forward/backward solves with the factors $\Lh$ and $\Uh$ that are
stored in adaptive precision. It is therefore important to analyze the accuracy of these solves.

\subsection{Algorithm and analysis}

A first option is to perform the solves entirely in high precision by converting the lower precision entries of $\Lh$ and $\Uh$ back
to high precision. This conversion can be implemented efficiently on-the-fly using
so-called memory accessors~\cite{ajlm26,gaq21,gjmm24a}. Since triangular solves are memory-bound
operations, the storage reduction achieved by the use of adaptive precision should translate, at least partially,
into execution time reductions, even if the solve operations are performed in high precision arithmetic.
Standard rounding error analysis of triangular
solves~\cite[Chap.~8]{high:ASNA2} shows that the backward error for the solves is then of order $\zeta$,
the unit roundoff of the compute precision.

Alternatively, it is of interest to consider whether the triangular solves can be performed in adaptive precision to reduce the
number of required conversions. Recall that, given the lower triangular system $\Lh y=b$ and the upper triangular system $\Uh x=y$,
the standard forward and backward substitution formulas give
\begin{align*}
y_i &= b_i - \sum_{j=1}^{i-1} \lh_{ij} y_j, \quad i=1\colon n,  \\
x_i &= \bigg(y_i - \sum_{j=i+1}^{n} \uh_{ij} x_j\bigg)/\uh_{ii}, \quad i=n\colon -1\colon 1.
\end{align*}
In our context, the factors can be expressed as the sum of $p$ disjoint sparse matrices of different precisions,
\begin{align*}
    \Lh &=  \Lh^{(1)} + \Lh^{(2)} + \dots + \Lh^{(p)}, \\
    \Uh &=  \Uh^{(1)} + \Uh^{(2)} + \dots + \Uh^{(p)},
\end{align*}
where $\Lh^{(i)}$ (respectively, $\Uh^{(i)}$) holds the subset of the entries of $\Lh$ (respectively, $\Uh$)
that are stored in precision $\tau_{i}$.
Algorithm~\ref{alg.solve} outlines evaluating these formulas in adaptive precision: essentially, terms $\lh_{ij} y_j$ or $\uh_{ij} x_j$ corresponding
to entries $\lh_{ij}^{(k)}$ or $\uh_{ij}^{(k)}$ in the same precision $\tau_k$ are summed  in precision $\tau_k$. Their partial sums are then converted and summed
in high precision $\tau_1$.
This reduces the number of conversions to at most $p$ per row (ignoring the conversions of the right-hand side and solution vectors).

\begin{algorithm}
\caption{Adaptive precision forward and backward substitution.\label{alg.solve}}
\textbf{Input:} \begin{tabular}[t]{l} $\Lh = \Lh^{(1)} + \dots + \Lh^{(p)} \in \R^{n\times n}$ in adaptive precision, \\
                $\Uh = \Uh^{(1)} + \dots + \Uh^{(p)} \in \R^{n\times n}$ in adaptive precision, \\
                $b\in\R^n$ in precision $\tau_1$ \end{tabular}\\
\textbf{Output:} $x\in\R^n$ in precision $\tau_1$ such that $\Lh\Uh x \approx b$

\begin{algorithmic}[1]
\FOR{$i=1\colon n$}
\FOR{$k=1\colon p$}
\STATE Initialize $y_i^{(k)} = 0$ 
\FOR{\textbf{each} $j<i$ s.t. $\lh^{(k)}_{ij} \ne 0$}
\STATE $y_i^{(k)} = y_i^{(k)} + \lh^{(k)}_{ij} y_j$ in precision $\tau_k$
\ENDFOR
\ENDFOR
\STATE $y_i = b_i - \sum_{k=1}^p y_i^{(k)}$ in precision $\tau_1$
\ENDFOR
\FOR{$i=n\colon -1\colon 1$}
\FOR{$k=1\colon p$}
\STATE Initialize $x_i^{(k)} = 0$
\FOR{\textbf{each} $j>i$ s.t. $\uh^{(k)}_{ij} \ne 0$}
\STATE $x_i^{(k)} = x_i^{(k)} + \uh^{(k)}_{ij} x_j$ in precision $\tau_k$
\ENDFOR
\ENDFOR
\STATE $x_i = \left( y_i - \sum_{k=1}^p x^{(k)}_i\right)/\uh_{ii}$ in precision $\tau_1$
\ENDFOR
\end{algorithmic}
\end{algorithm}

Unfortunately, the next result shows that Algorithm~\ref{alg.solve} is not necessarily stable.

\begin{theorem}\label{thm.solve}
Let Algorithm~\ref{alg.solve} be used to solve the linear system $\Lh \Uh x =b$,
where the $\Lh$ and $\Uh$ factors have been obtained by Algorithm~\ref{alg.ILUT}.
If the algorithms run to completion, the computed solution $\xh$ satisfies
\begin{equation}
\begin{aligned}\label{eq.bnd-solve}
&(\Lh+\Delta L)(\Uh+\Delta U)\xh = b, \\
&\quad \abs{\Delta L} \le \gamma_{\ceps+p}^{(1)}\abs{\Lh} + \cepsp\,\epsabs J_L\abs{\diag(\Uh)^{-1}}, \\
&\quad \abs{\Delta U} \le \gamma_{\ceps+p}^{(1)}\abs{\Uh} + \cepsp\, \epsabs J_U,
\end{aligned}
\end{equation}
where $J_L$ and $J_U$ are matrices of ones and zeros with the same sparsity pattern as $\Lh$ and $\Uh$, respectively, 
$\ceps \le n$ is the maximum number of nonzeros on any row or column of $\Lh$ and $\Uh$,
$\gamma_{\ceps+p}^{(1)} = (\ceps+p)\tau_1/(1-(\ceps+p)\tau_1)$,
and $\cepsp = (\ceps+p+1)/(1-(\ceps+p+1)\tau_p)$.
\end{theorem}

\begin{proof}
Consider first the forward substitution $\Lh y =b$.
The computed $\yh_i^{(k)}$ satisfies
\[
\yh_i^{(k)} = \sum_{j=1}^{i-1} \lh_{ij}^{(k)} \yh_j (1+\theta_{ij}^{(k)}), \quad \abs{\theta_{ij}^{(k)}} \le \gamma_{\ceps}^{(k)},
\]
where we define $\gamma_{\ceps}^{(k)} = \ceps\tau_k/(1-\ceps\tau_k)$.
The computed $\yh_i$ thus satisfies
\begin{align*}
\yh_i(1+\xi_i^{(p+1)}) &= b_i - \sum_{k=1}^p \yh_i^{(k)}(1+\xi_i^{(k)}), \quad\abs{\xi_i^{(k)}} \le \gamma_p^{(1)} \\
&=b_i - \sum_{j=1}^{i-1}\sum_{k=1}^p \lh_{ij}^{(k)} \yh_j (1+\theta_{ij}^{(k)})(1+\xi_i^{(k)}).
\end{align*}
Note that the $\Lh^{(k)}$ have disjoint sparsity patterns and so
$\sum_{k=1}^p \lh_{ij}^{(k)} = \lh_{ij} = \lh_{ij}^{(k)}$ for some $k$. Moreover,
since $\lh_{ii} = 1$, we have
\begin{align*}
b_i &= \sum_{j=1}^{i} \lh_{ij}^{(k)} \yh_j (1+\theta_{ij}^{(k)})(1+\xi_i^{(k)}) \\
    &=  \sum_{j=1}^{i} \lh_{ij} \yh_j + \sum_{j=1}^{i} \lh_{ij}^{(k)} \yh_j (\theta_{ij}^{(k)} + \xi_i^{(k)} + \theta_{ij}^{(k)}\xi_i^{(k)}).
\end{align*}
In matrix form, we obtain
\[
 b = (\Lh + \Delta L)\yh, \quad \abs{\Delta \ell_{ij}}  
 \le (\gamma_{\ceps}^{(k)} + \gamma_p^{(1)} + \gamma_{\ceps}^{(k)}\gamma_p^{(1)})  \abs{\lh_{ij}^{(k)}}
 \le \gamma_{\ceps+p}^{(k)}                                                           \abs{\lh_{ij}^{(k)}},
\]
where we have used~\cite[Lemma~3.3]{high:ASNA2}
\[
\gamma_{\ceps}^{(k)} + \gamma_p^{(1)} + \gamma_{\ceps}^{(k)}\gamma_p^{(1)} 
\le
\gamma_{\ceps}^{(k)} + \gamma_p^{(k)} + \gamma_{\ceps}^{(k)}\gamma_p^{(k)} 
\le
\gamma_{\ceps+p}^{(k)}.
\]
If $k\ge 2$, $\abs{\ell_{ij}^{(k)}\uh_{jj}} \le \epsabs/\tau_k$, so
$\abs{\lh_{ij}^{(k)}\uh_{jj}} \le (1+\tau_k)\epsabs/\tau_k$ and 
\[
 b = (\Lh + \Delta L)\yh, \quad \abs{\Delta \ell_{ij}} 
\le \gamma_{\ceps+p}^{(1)}\abs{\lh_{ij}} + \frac{(1+\tau_k)\epsabs}{\tau_k \abs{\uh_{jj}}} \gamma_{\ceps+p}^{(k)}
\le \gamma_{\ceps+p}^{(1)}\abs{\lh_{ij}} + \frac{\cepsp\epsabs}{\abs{\uh_{jj}}} 
\]
with $\cepsp = (\ceps+p+1)/(1-(\ceps+p+1)\tau_p)$, 
where we have used $(1+\tau_k)\gamma_{\ceps+p}^{(k)} \le \gamma_{\ceps+p+1}^{(k)}$.
In matrix form, this finally yields
\begin{equation}\label{eq.DeltaL}
\abs{\Delta L} \le \gamma_{\ceps+p}^{(1)}\abs{\Lh} + \cepsp\,\epsabs J_L\abs{\diag(\Uh)^{-1}},
\end{equation}
where $J_L$ is the matrix with ones in the same positions as the nonzeros of $\Lh$ and zeros elsewhere.

The analysis for the backward substitution is similar and yields the bound 
\begin{equation}\label{eq.DeltaU}
(\Uh+\Delta U)\xh = \yh, \quad \abs{\Delta U} \le \gamma_{\ceps+p}^{(1)}\abs{\Uh} + \cepsp\, \epsabs J_U.
\end{equation}
Overall, the computed solution $\xh$ therefore satisfies
\begin{equation}\label{eq.DeltaLU}
(\Lh+\Delta L)(\Uh+\Delta U)\xh = b,
\end{equation}
with bounds on $\Delta L$ and $\Delta U$ given in \eqref{eq.DeltaL} and \eqref{eq.DeltaU}.
\end{proof}

This is not an entirely satisfying bound, since we would like the perturbations $\Delta L$ and $\Delta U$ to be relative
to $\abs{\Lh}$ and $\abs{\Uh}$ as for a standard triangular solve, but this is not the case because of the terms
proportional to $\epsabs$.
In fact, if we rewrite \eqref{eq.bnd-solve} as a bound on the residual $\Lh\Uh\xh-b$, we obtain
to first order
\[
\abs{\Lh\Uh\xh -b} \lesssim (\ceps+p)\left( 2\tau_1\abs{\Lh}\abs{\Uh}
+ \epsabs \abs{\Lh}J_U + \epsabs J_L\abs{\diag(\Uh)^{-1}\Uh} \right)\abs{\xh} 
\]
Taking norms and using $\norminf{J_L} \le \ceps$ and $\norminf{J_U} \le \ceps$ yields 
\begin{equation}\label{eq.bnd.x}
\frac{\norminf{\Lh\Uh\xh-b}}{\norminf{\xh}} \lesssim 
(\ceps+p)\left(2\tau_1\norminf{\abs{\Lh}\abs{\Uh}} + \ceps\epsabs\norminf{\Lh} + \ceps\epsabs\norminf{\diag(\Uh)^{-1}\Uh} \right).
\end{equation}

As mentioned in Remark~\ref{rem.epschoice}, in practice a typical choice for $\epsabs$
would be $\epsabs = \epsrel\norminf{A}$---which yields a normwise relative backward error
of order $\epsrel$ for the ILU factorization (Algorithm~\ref{alg.ILUT}). Plugging this value in \eqref{eq.bnd.x}
yields a normwise relative backward error bound for the solution of $\Lh\Uh x=b$:
\begin{equation}\label{eq.bnd.x2}
\frac{\norminf{\Lh\Uh\xh - b}}{\norminf{A}\norminf{\xh}}  \lesssim
(\ceps+p)\left(2\tau_1\frac{\norminf{\abs{\Lh}\abs{\Uh}}}{\norminf{A}} + \ceps\epsrel\norminf{\Lh} + \ceps\epsrel \norminf{\diag(\Uh)^{-1}\Uh} \right).
\end{equation}

This bound can only be considered satisfactory if the terms in front of $\tau_1$ and $\epsrel$ are small.
We can guarantee the following:
\begin{itemize}
\item $\norminf{\abs{\Lh}\abs{\Uh}}/\norminf{A}$ is bounded by $\ceps^2\xi^{-1}\rho_n$
if threshold partial pivoting (with threshold $\xi$) is used, 
where the growth factor $\rho_n$ is often small;
\item $\norminf{\Lh}$ is bounded by $\ceps\xi^{-1}$ if threshold partial pivoting is used;
\item $\norminf{\diag(\Uh)^{-1}\Uh}$ is bounded by $\ceps\xi^{-1}$ if either threshold rook or complete pivoting 
is used, but can be arbitrarily large if only partial pivoting is used.
\end{itemize}
Partial pivoting selects the largest entry in the current column (or an entry within a given threshold $\xi$ of the largest entry);
rook pivoting iterates to find an entry that is the largest in both its own row and its own column; see~\cite[Chap.~9]{high:ASNA2}.
While (threshold) partial pivoting is commonly employed, (threshold) rook pivoting is more costly and less common.
Therefore, the most problematic term in bound \eqref{eq.bnd.x2} is $\norminf{\diag(\Uh)^{-1}\Uh}$.

\subsection{Numerical example}
\label{sec:solves-example}

We have tested Algorithm~\ref{alg.solve} on a wide range of SuiteSparse matrices and found it to be
much more accurate than \eqref{eq.bnd.x2} suggests on most of them. However, for
a few examples it does cause instability unless rook pivoting is used, which confirms that the bound can be descriptive in some cases.
We report one such example in Table~\ref{tab:solve}. 

\begin{table}[htbp]
    \centering
    \caption{Accuracy of Algorithm~\ref{alg.solve} on matrix nnc666 with $\epsrel=10^{-16}$.}
    \label{tab:solve}
    \footnotesize
    \begin{tabular}{lcccc|ccc}
    \toprule
    Pivot & $\ceps$ & $\frac{\norminf{\abs{\Lh}\abs{\Uh}}}{\norminf{A}}$ & $\norminf{\Lh}$ & $\norminf{\diag(\Uh)^{-1}\Uh}$ & $\frac{\norminf{\Lh\Uh\xh - b}}{\norminf{A}\norminf{\xh}}$ & RHS of \eqref{eq.bnd.x2} & $\frac{\text{RHS of } \eqref{eq.bnd.x2}}{\ceps^2}$ \\
    \midrule
    Partial & $5{\times}10^2$ & $1{\times}10^1$ & $3{\times}10^1$ & $1{\times}10^7$ & $3{\times}10^{-10}$ & $3{\times}10^{-4}$ & $1{\times}10^{-9\phantom{0}}$  \\
    Rook    & $5{\times}10^2$ & $2{\times}10^1$ & $3{\times}10^1$ & $3{\times}10^1$ & $1{\times}10^{-14}$ & $2{\times}10^{-9}$ & $6{\times}10^{-15}$ \\
    \bottomrule
    \end{tabular}
\end{table}

While $\frac{\norminf{\abs{\Lh}\abs{\Uh}}}{\norminf{A}}$ and $\norminf{\Lh}$ are bounded by small constants for both partial and rook pivoting, $\norminf{\diag(\Uh)^{-1}\Uh}$ is much larger when partial pivoting is used. This leads to a loss of accuracy of more than four digits with partial pivoting, which is reflected in the bound
\eqref{eq.bnd.x2}. The bound itself is several orders of magnitude larger than the actual error, but this is mainly due to the dimensional constant $\ceps^2$ ($\approx 2\times 10^5$ for this matrix), which is known to be pessimistic~\cite{hima19a,begm}; in fact, dropping $\ceps^2$ from the bound (last column of Table~\ref{tab:solve}) makes it of the same order as the measured error.

We conclude from this experiment that the bound \eqref{eq.bnd.x2} can be descriptive of real instabilities when Algorithm~\ref{alg.solve} is used without rook pivoting. We thus recommend keeping the triangular solves
in higher precision, except in an environment where the required conversions become more costly than the use of (threshold) rook pivoting.

\begin{remark} \label{rem.tri_err}
    We briefly mention another more sophisticated approach to reduce the number of high-precision conversions while avoiding the need for rook pivoting.
    Inspecting the analysis of Algorithm~\ref{alg.solve} shows that instabilities are created by accumulating in low precision entries $\lh_{ij}$ belonging to a column with a small diagonal
    pivot $\uh_{jj}$. Modifying the algorithm to convert to high precision entries $\lh_{ij}$ such that $\abs{\lh_{ij}}\tau_k > \epsabs$ removes the 
    $\norminf{\diag(\Uh)^{-1}\Uh}$ from the error bound and cures the instability.
    On matrix nnc666 of Table~\ref{tab:solve}, for example, this modification restores a backward error of order $4\times 10^{-13}$ when using partial pivoting, at the price of additional conversions.
\end{remark}

\section{Adaptive precision threshold-based incomplete Cholesky}
\label{sec:IC}

\begin{algorithm}
\caption{Adaptive precision IC (\APIC)\label{alg.IC}}
\textbf{Input:}  \begin{tabular}[t]{l}symmetric, positive-definite matrix $A\in\R^{n\times n}$, target accuracy $\epsabs$, \\
                 storage precisions $\{\tau_1,\ldots,\tau_p\}$, compute precision $\zeta$  \end{tabular} \\
\textbf{Output:} $\Lh\Lh^{T}\approx A$

\begin{algorithmic}[1]
\FOR{$k=1\colon n$}
\STATE $\ell_{kk} = \sqrt{a_{kk} - \sum_{i=1}^{k-1} \lh_{ki}^2}$ in precision $\zeta$ \label{line.IC.lkk}
\STATE $\tau = \max\left\{\tau \in \{\tau_1,\ldots,\tau_p\} \colon \abs{\ell_{kk}^2} \le \epsabs/(2\tau + \tau^2) \right\}$
\STATE $\lh_{kk} = \fl_\tau(\ell_{kk})$ \label{line.IC.hlkk}
\FOR{$i=k+1\colon n$}
\STATE $\ell_{ik} = \big(a_{ik} - \sum_{j=1}^{k-1} \lh_{ij}\lh_{kj}\big)/\lh_{kk}$ in precision $\zeta$
\STATE $\tau = \max\left\{\tau \in \{\tau_1,\ldots,\tau_p\} \colon \abs{\ell_{ik}\lh_{kk}} \le \epsabs/\tau \right\}$
\STATE $\lh_{ik} = \fl_\tau(\ell_{ik})$
\ENDFOR
\ENDFOR
\end{algorithmic}
\end{algorithm}

So far we have considered unsymmetric systems and general ILU factorizations. We now focus on SPD systems, for which we use an incomplete Cholesky (IC) factorization.

Algorithm~\ref{alg.IC} describes an adaptive precision variant of threshold-based IC (\APIC), where the entries retained in the factor are stored in a precision inversely proportional to their magnitude so that the backward error remains bounded by a given tolerance.
The algorithm is analogous to Algorithm~\ref{alg.ILUT}, but with one noteworthy difference:
since the diagonal entries are stored as a square root, the criterion to choose their precision is slightly more stringent (replacing the unit roundoff $\tau$ by $2\tau + \tau^2$).

The following theorem proves the error result for the proposed \APIC and is analogous to Theorem~\ref{thm.ILUT}.

\begin{theorem}\label{thm.IC}
    If Algorithm~\ref{alg.IC} applied to the SPD matrix $A\in\R^{n\times n}$ runs to completion, the computed Cholesky factor satisfies
    \begin{equation}
        A+E = \Lh\Lh^{T},\quad \abs{E} \le (1 + \gamma_{\ceps + 1})\epsabs J + \gamma_{\ceps + 1} \abs{\Lh}\abs{\Lh}^{T},
    \end{equation}
    where $J$ is the matrix of ones,
    $\ceps \le n$ is the maximum number of nonzeros on any row of $\Lh$,
    and $\gamma_{\ceps+1} = (\ceps+1)\zeta / \big(1 - (\ceps+1) \zeta\big)$.
\end{theorem}

\begin{proof}
    Let $\ell'_{kk} = (a_{kk} - \sum_{i=1}^{k-1} \lh^2_{ki})^{1/2}$.
    By \cite[Problem~10.3]{high:ASNA2}, the computed $\lh_{kk}'$ satisfies
    \[
        (\lh'_{kk})^2(1 + \theta^{(k)}_{c_{\epsrel}+1}) = a_{kk} - \sum_{i=1}^{k-1} \lh^2_{ki}(1 + \theta^{(i)}_{c_{\epsrel}}).
    \]
    Let $\lh_{kk} = \fl_\tau(\lh_{kk}') = \lh_{kk}' + e_{kk}$, with $\abs{e_{kk}} \le \tau\abs{\lh'_{kk}}$.
    We have
    \[
    (\lh_{kk})^2 = (\lh_{kk}' + e_{kk})^2 = (\lh_{kk}')^2 + e_{kk}^2 + 2\lh_{kk}'e_{kk} = (\lh_{kk}')^2 + f_{kk},
    \]
    with $\abs{f_{kk}} \le (2\tau+\tau^2)(\lh'_{kk})^2 \le \epsabs$, by construction.
    Therefore
    \[
        (\lh_{kk}^2-f_{kk})(1 + \theta^{(k)}_{c_{\epsrel}+1}) = a_{kk} - \sum_{i=1}^{k-1} \lh^2_{ki}(1 + \theta^{(i)}_{c_{\epsrel}}) 
    \]
    yields
    \[
        a_{kk} = \sum_{i=1}^{k} \lh^2_{ki}(1 + \theta^{(i)}_{c_{\epsrel}+1}) - f_{kk}(1 + \theta^{(k)}_{c_{\epsrel}+1}).
    \]
    
 For the off-diagonal entries, $a_{ik}, 1 \leq k < i \leq n$, we can similarly prove that
  \[
  \sum_{j=1}^{k}\lh_{ij} \lh_{kj} = a_{ik} + \delta a_{ik}, \quad \abs{\delta a_{ik}} \le (1 + \gamma_{\ceps+1})\epsabs + \gamma_{\ceps+1} \sum_{j=1}^{k} \abs{\lh_{ij}} \abs{\lh_{kj}}.
  \]
  The proof is completed by writing the errors in matrix form.
\end{proof}

\section{Adaptive precision memory-limited IC factorizations}
\label{sec:memlim}

So far, we have considered threshold-based factorizations, which provide a direct control over the backward error through the threshold parameter $\epsrel$, as
in Theorems~\ref{thm.ILUT} and \ref{thm.IC}. However, they do not provide an a priori control over the storage required by the factors.

To control the storage more directly, we may consider
structure-based methods through the use of an initial symbolic factorization phase 
to find the location of permissible entries using only 
the sparsity pattern ${\mathcal S}_A$ of $A$.
A widely used example is ILU(0) 
(or IC(0) for SPD systems)
that allows the incomplete factors to
have entries only in the positions corresponding to ${\mathcal S}_A$.
The success of structure-based methods is highly problem dependent and, in general,
they do not guarantee any control over the error since arbitrarily large entries may be dropped.
Memory-based methods can provide a more robust alternative, in which the memory available for
the incomplete factorization is prescribed and, within each column, only the largest entries are retained.
While they do not provide a direct control over the error like threshold-based methods, the error can still
be indirectly controlled, since increasing the available memory allows smaller entries to be retained.
A discussion of the various classes of incomplete factorization preconditioners (with references) is given in \cite{sctu:2023}.

In this section, we use adaptive precision storage within memory-limited incomplete factorizations. 
Our experiments are based on the \texttt{HSL\_MI28} package, which computes a memory-limited IC factorization for sparse SPD systems.

\subsection{Memory-limited factorizations}

We employ a memory-limited IC
factorization approach that uses additional memory during the construction of
$\Lh$ that is subsequently discarded at the end of the factorization. The aim is to obtain a high
quality preconditioner while maintaining sparsity and allowing the
user to control the memory usage.
Consider the decomposition 
\begin{equation}\label{eq:memory_limited}
 A = (\Lh+ \Rh) \, (\Lh+ \Rh)^T - E,
\end{equation}
where $\Lh$ is  lower triangular with positive diagonal
entries, $\Rh$ is strictly lower
triangular and $ E = \Rh \Rh^T$
is the error matrix.
At step $j-1$ ($j \ge 2$) of the factorization, the first column of
the Schur complement is split into the sum\footnote{Here we use the standard section notation, that is, $L_{i:k,j}$ denotes the entries in rows $i$ to $k$ of column $j$.}
$\Lh_{j:n,j}+\Rh_{j:n,j},$
where $\Lh_{j:n,j}$
contains the entries that are retained in  column $j$
of $\Lh$, 
the diagonal entry $\Rh_{j,j}$ is zero, and $\Rh_{j+1:n,j}$
contains the entries that are  not included in $\Lh$.
In a complete factorization,
the Schur complement is updated by subtracting
$$(\Lh_{j+1:n,j}+\Rh_{j+1:n,j}) \, (\Lh_{j+1:n,j}+\Rh_{j+1:n,j})^T .$$
However, the incomplete factorization discards the term
\begin{equation*}
E^{(j)} =  \Rh_{j+1:n,j} \, \Rh_{j+1:n,j}^T.
\end{equation*}
Thus, $E^{(j)}$ is implicitly added to $A$ and because $E^{(j)}$
is positive semidefinite, in exact arithmetic, the approach does not break down.
An obvious choice is for the largest entries in the column to be retained in $\Lh$.
More fill entries are used in constructing $\Lh$ than in a standard factorization, allowing
the structure of the complete factorization to be followed more closely.

The computational complexity
and memory needed can be reduced by limiting how many entries are allowed in each
column of $\Lh$ and $\Rh$. In the package {\tt HSL\_MI28},  user-chosen parameters determine the memory needed as well as the work involved in the factorization:  
{\tt lsize} controls the number of fill entries in each column of $\Lh$ 
while ${\tt rsize}$ controls the number of entries in each column of $\Rh$. {\tt HSL\_MI28} also incorporates the
use of drop tolerances
$\delta_1 \ge 0$ and $\delta_2 \ge 0$ for the entries of $\Lh$ and $\Rh$, respectively.
Discarding entries can cause the IC factorization to break down because a diagonal entry 
(the pivot) becomes too small (or negative).
In this case,  a (small)  shift $\alpha >0$ is chosen and the
IC factorization of the symmetrically scaled and shifted matrix $\Theta^{-1}A\Theta^{-1} + \alpha I$ is computed, 
where $\Theta$ is a diagonal scaling matrix. {\tt HSL\_MI28}
handles the use of a shift automatically and includes options for 
preordering and prescaling $A$ \cite{sctu:2014a,sctu:2014b}.

Let us now discuss how to adapt Algorithm~\ref{alg.IC} to this memory-limited case.
An important observation is that even though memory-limited factorizations do not directly bound the size of the dropped entries,
their accuracy is still determined by it. Thus, Theorem~\ref{thm.IC} still applies to memory-limited factorizations,
where $\epsabs$ is no longer a user-controlled parameter but rather the effective largest magnitude $\abs{\ell_{ik}\lh_{kk}}$ for any dropped entry $\ell_{ik}$,
which is indirectly controlled by the \texttt{lsize} parameter.

The key consequence of this observation is that the adaptive precision memory-limited factorization should be constructed
using the same criterion as in Algorithm~\ref{alg.IC}, which by design balances the introduced rounding errors with the errors
from dropping entries of magnitude at most $\epsabs$. One caveat is that $\epsabs$ is only known a posteriori; in principle, we could compute
the factors in uniform precision to estimate $\epsabs$ and then switch the factors (or recompute them) in adaptive precision. In practice, though,
it is sufficient to have a rough estimate of the order of magnitude of $\epsabs$ in order to balance approximately the two sources of error, as we will illustrate in the next section.

\begin{remark}\label{rem.dynamic_epsa}
    A more sophisticated approach that only requires one factorization would be to dynamically update $\epsabs$
    as a running maximum of the magnitude of the largest entry dropped so far. One could hope that $\epsabs$
    would be quickly estimated from the elimination of the first few columns. Investigating this idea is left for future work.
\end{remark}

\begin{remark}
    Our error analysis applies to memory-limited factorizations with $\texttt{rsize}=0$. 
    Even though, in practice, the $\Rh$ factor usually helps reduce the number of iterations,
    quantifying theoretically its effect on the backward error is outside our scope.
\end{remark}

\subsection{Numerical example}
\label{sec:memlim-example}

Figure~\ref{fig:memlim} illustrates how the performance of the preconditioner varies as both the fill
parameter \texttt{lsize} and the adaptive tolerance $\epsrel$ (where $\epsabs = \epsrel\norm{A}$) are varied
for the matrix \texttt{bcsstk13}.
Throughout this section we choose $\texttt{rsize}=\texttt{lsize}$ when generating the incomplete factors.
Each cell of the heatmap corresponds to a memory-limited IC factorization
computed by a modified version of {\tt HSL\_MI28} in which entries are adaptively stored in lower precisions.
The color in Figure~\ref{fig:memlim_error} shows an estimate of the relative backward error of
the computed factorization,
$\norm{A - \Lh\Lh^T}_F/\norm{A}_F$,
obtained efficiently via the randomized Hutchinson trace estimator \cite{hutchinson_01_1989,martinsson_tropp_05_2020}.
Figure~\ref{fig:memlim_iters} shows the corresponding number of PCG iterations.

\begin{figure}
    \centering
    \begin{subfigure}[c]{0.49\textwidth}
        \includegraphics[width=\linewidth]{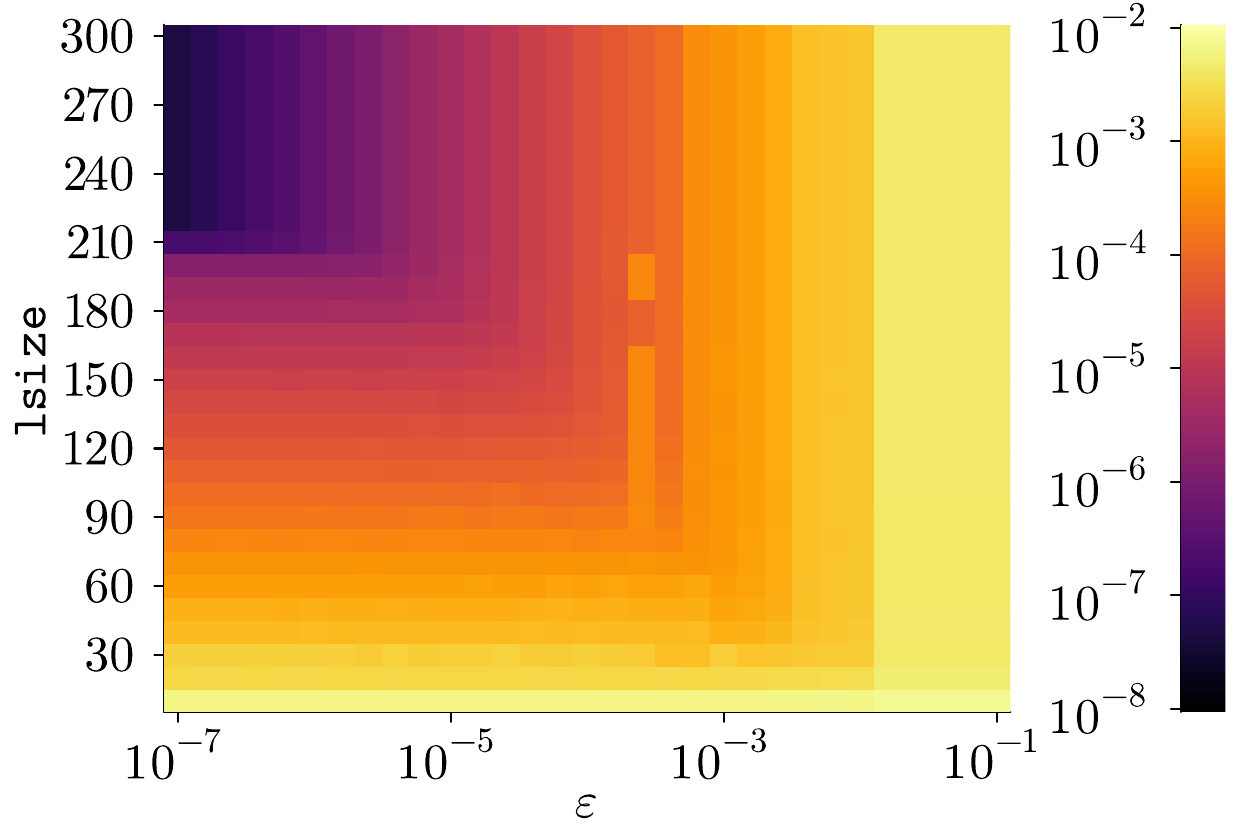}
        \caption{Backward error}
        \label{fig:memlim_error}
    \end{subfigure}
    \hfill
    \begin{subfigure}[c]{0.49\textwidth}
        \includegraphics[width=\linewidth]{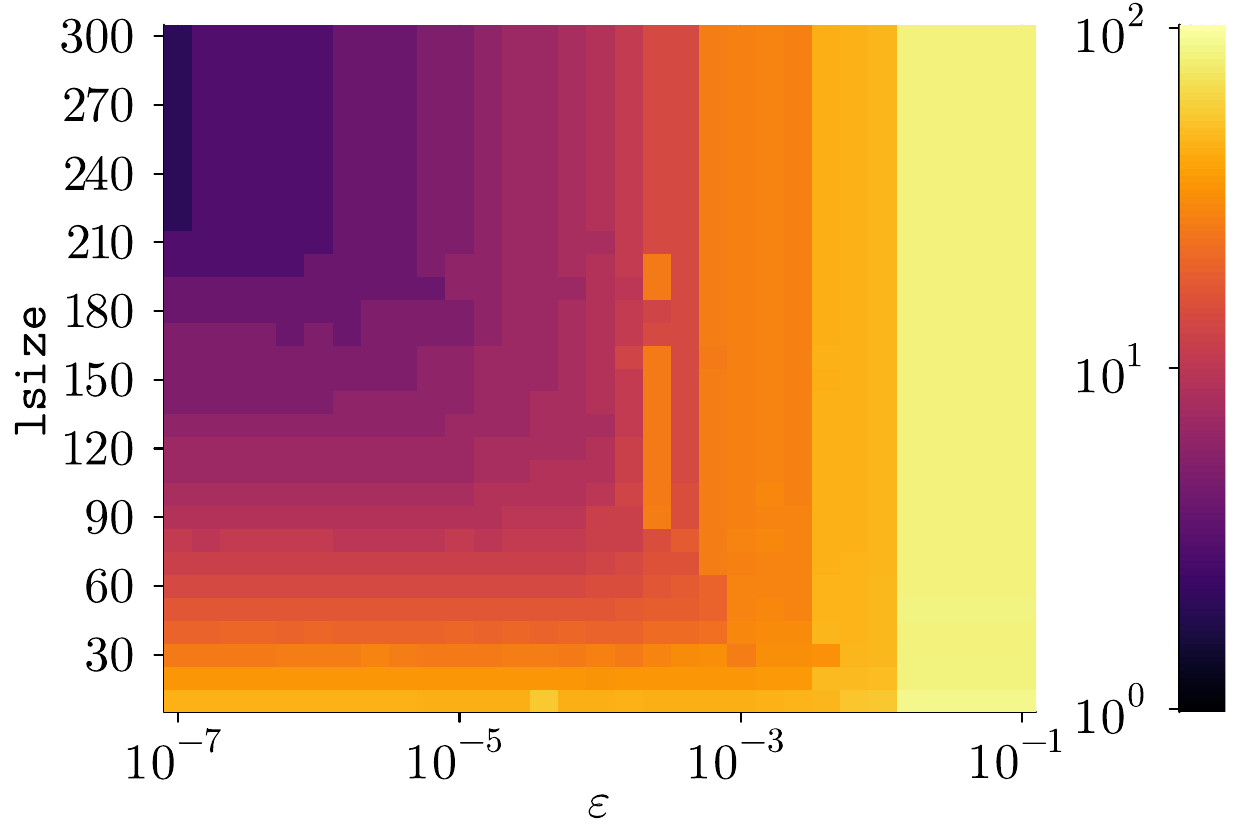}
        \caption{PCG iterations}
        \label{fig:memlim_iters}
    \end{subfigure}
    \caption{Results for the memory-limited adaptive precision preconditioner applied to matrix \texttt{bcsstk13}.
             Each cell corresponds to an IC factorization with the given fill parameter \texttt{lsize}
             and adaptive tolerance $\epsrel$, used as a preconditioner for PCG.}
    \label{fig:memlim}
\end{figure}

These plots highlight the balance between the two sources of error: the rounding error introduced
by storing entries in lower precision (controlled by $\epsrel$), and the memory-limited error from
restricting the number of entries per column of $\Lh$ (controlled by \texttt{lsize}).
When $\epsrel$ is large, increasing \texttt{lsize} does not improve the preconditioner because the
rounding error dominates.
Conversely, when \texttt{lsize} is small, the memory-limited error dominates and the iteration count is
insensitive to $\epsrel$ over a wide range, confirming that storing the factors in adaptive precision
does not degrade performance in this memory-limited regime.

\begin{table}[tbp]
\centering
\caption{Storage for \texttt{bcsstk13}: uniform fp64, fp32 and rp16 versus adaptive precision (\APIC) with $\epsrel = 10^{-4}$. ``FP'' is the storage for the floating-point values only, ``Int'' is the storage for the integer 4-byte row indices and column pointers assuming a CSC sparse matrix format. ``\#its'' is the PCG iteration count.}
\label{tab:memlim_storage_1e-4}
\begin{tabular}{lrrrrr}
\toprule
Method & \texttt{lsize} & \multicolumn{3}{c}{size($\Lh$) in kB} & \#its \\
& & FP & Int & Total \\
\midrule
Uniform fp64 & 40 & 951 & 483 & 1434 & 21 \\
Uniform fp32 & 40 & 475 & 483 & 959 & 21 \\
Uniform rp16 & 40 & 238 & 483 & 721 & 83 \\
\APIC & 40 & 257 & 483 & 740 & 20 \\
\APIC & 60 & 332 & 634 & 966 & 15 \\
\APIC & 100 & 480 & 929 & 1409 & 10 \\
\bottomrule
\end{tabular}
\end{table}

Table~\ref{tab:memlim_storage_1e-4} seeks to put these results into context.
In the experiments using {\tt HSL\_MI28} (in uniform fp64) reported in \cite{sctu:2014a},  $\texttt{lsize}$ is chosen to be at most 40.
At this fill level, for matrix \texttt{bcsstk13} 
\APIC with a tolerance of $\epsrel = 10^{-4}$ can be used without degrading the PCG iteration count while reducing the storage requirement for the floating-point values by $73\%$ compared to fp64; this translates to an overall $48\%$ storage reduction when also accounting for the integer overhead of the sparse format. Conversely, another point of view is that thanks to the use of \APIC we could increase
\texttt{lsize} to $100$ and match the fp64 storage, while halving the iteration count.

To be fair, we should also compare \APIC to uniform reduced precision variants. Uniform fp32 performs well on this matrix, as it achieves the same iteration count as fp64; however, \APIC performs even better, further reducing the floating-point storage by 46\%, which translates to a 23\% total reduction. Alternatively, taking $\texttt{lsize}=60$ matches the fp32 storage and reduces the iteration count from 21 to 15. Uniform rp16 performs poorly: it achieves a storage that is only marginally smaller than \APIC, while increasing the iteration count by $4\times$.

Finally, these results illustrate that the relative fraction of storage required by the integers increases when using adaptive precision, and can be quite significant, especially for larger values of \texttt{lsize}. This motivates investigating more compact sparse formats with compressed representations of the indices as future work.

\subsection{Numerical experiments on a large test set}
\label{sec:exps}

We now assess the memory-limited \APIC approach across a broader test set.
We use the 153 symmetric positive-definite matrices from the SuiteSparse Matrix Collection
that were employed in the original {\tt HSL\_MI28} paper \cite{sctu:2014a}.
For each matrix, we run our modified {\tt HSL\_MI28} code with fill parameters $\texttt{lsize} \in \{20, 60, 100, 140, 180, 220\}$
and $\texttt{rsize} = \texttt{lsize}$,
drop tolerances $\delta_1 = \delta_2 = 0$, and adaptive tolerance
$\epsrel \in \{10^{-2}, 10^{-4}, 10^{-6}, 10^{-8}, 10^{-10}\}$,
alongside uniform fp64, fp32, and rp16.
We take $\texttt{rsize} = \texttt{lsize}$ because the factor $\Rh$ is discarded once $\Lh$ has been
computed, so it does not change the final storage or the per-iteration cost of the preconditioner, but
using $\texttt{rsize}>0$ generally improves the preconditioner quality.
In a test on a subset of $14$ matrices with $\texttt{lsize} \in \{20, 100\}$,
setting $\texttt{rsize} = 0$ increased the number of PCG iterations by a median of
$10$--$20\%$ and caused one additional solver failure (iteration count exceeded the limit).

\begin{figure}[ht]
    \centering
    \begin{subfigure}[c]{0.65\textwidth}
        \centering
        \includegraphics[width=\linewidth]{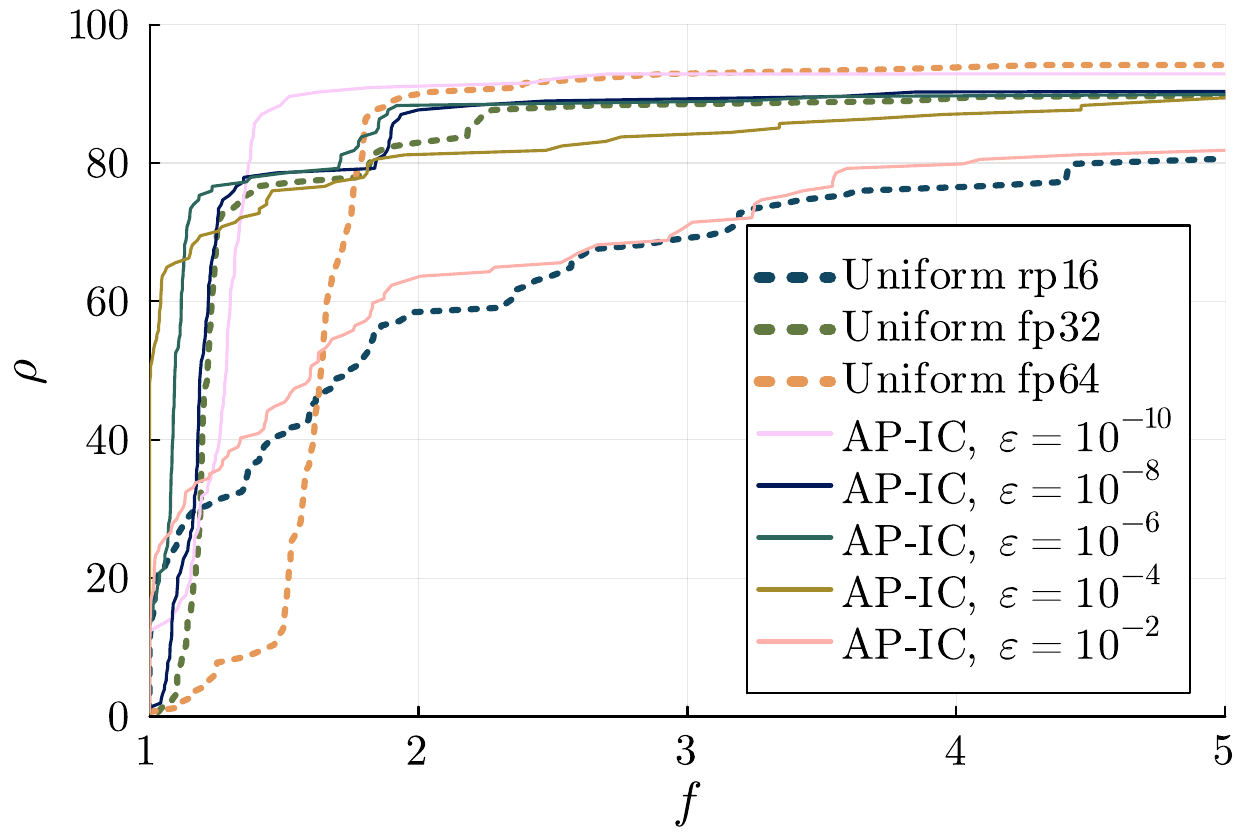}
        \caption{\texttt{lsize} $= 20$}
        \label{fig:memlim_perf_prof_20}
    \end{subfigure}

    \begin{subfigure}[c]{0.65\textwidth}
        \centering
        \includegraphics[width=\linewidth]{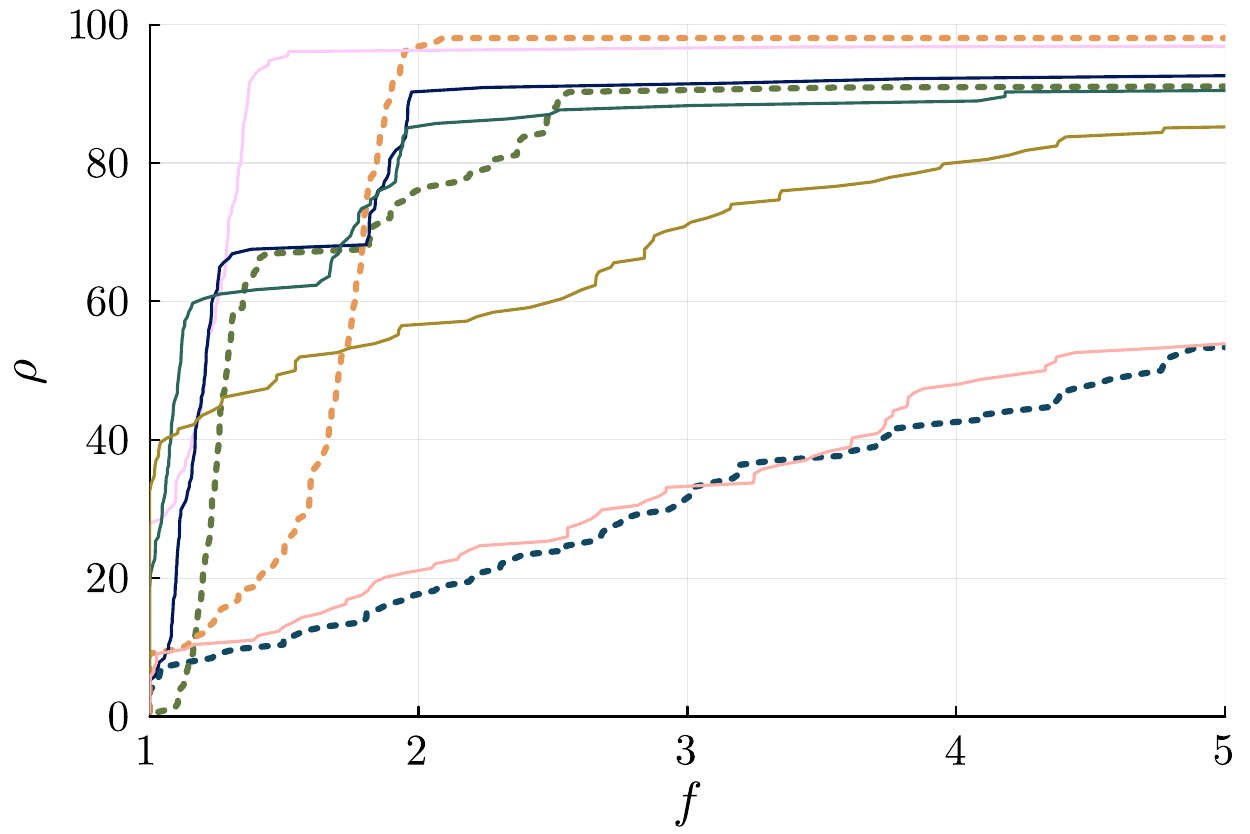}
        \caption{\texttt{lsize} $= 100$}
        \label{fig:memlim_perf_prof_100}
    \end{subfigure}
    \caption{Performance profiles comparing the cost metric \eqref{eqn:cost} for adaptive precision (\APIC) {\tt HSL\_MI28} with varying $\epsrel$ against uniform fp64, fp32, and rp16 baselines.
             For each problem, the cost of a method is normalized with respect to the lowest cost of any method in the same plot (same \texttt{lsize}); $\rho$ is the percentage of the test set for which the normalized cost is less than or equal to $f$.    \label{fig:memlim_perf_prof}
              }
\end{figure}

Figure~\ref{fig:memlim_perf_prof} presents performance profiles \cite{dolan_more_01_2002}
using the cost metric
\begin{equation}\label{eqn:cost}
    \text{cost} = \#\mathrm{its} \cdot (2 \times \text{size}(\Lh) + \text{size}(A)),
\end{equation}
where $\#\mathrm{its}$ is the PCG iteration count 
and $\text{size}(\Lh)$, $\text{size}(A)$ are the storage costs in bytes for the factor $\Lh$ and matrix $A$, respectively.
The storage costs assume a CSC sparse matrix format and include the integer overheads for row indices and column pointers.
This cost metric approximates the total floating-point work of the preconditioned solve.

We first analyze using $\texttt{lsize}=20$ (Figure~\ref{fig:memlim_perf_prof_20}). 
Uniform rp16 and \APIC with $\epsrel=10^{-2}$ are the most aggressive variants, which allows them to achieve the lowest cost for about 20\% of the problems; however, for the remaining problems, they require significantly more iterations and a higher cost to converge, and fail to converge to the requested tolerance for a larger proportion of the problems compared to the other variants.
Uniform fp64 and, to a lesser extent, uniform fp32 are the most conservative variants; they converge for most problems in the test set, but also require a higher cost on most of them: for example, the fp64  cost is at least 40\% higher than the best variant for about 70\% of the problems.
Overall, the \APIC variants achieve the best tradeoff between cost and robustness, especially for intermediate target accuracies such as $\epsrel = 10^{-4}$ and $\epsrel = 10^{-6}$, for which they achieve the lowest cost for a large proportion of the problems while still converging on almost as many problems as the most conservative variants.

Next, in the case $\texttt{lsize} = 100$ (Figure~\ref{fig:memlim_perf_prof_100}), the results lead to the same general conclusions, with a few noteworthy differences. First, the overall best performing \APIC variants use smaller target accuracies $\epsrel \le 10^{-6}$. This is consistent with the single-matrix analysis in section~\ref{sec:memlim-example}: at higher fill levels the memory-limited error is smaller, so a tighter $\epsrel$ is needed to keep the rounding error from dominating.
Second, the difference in robustness between conservative variants and aggressive ones (such as uniform rp16 or \APIC with a large $\epsrel$) is exacerbated; for example, uniform rp16 either fails to converge or costs more than five times as much as the best variant on nearly half the test set. This is because, with a  larger \texttt{lsize}, the conservative variants converge on the harder test problems,
while the aggressive ones remain limited by the rounding errors.

For the ten largest problems by number of nonzeros, we provide detailed results listed in Table~\ref{tab:memlim_largest}. We report the storage, PCG iterations, and cost \eqref{eqn:cost} (normalized relative to uniform fp64 with $\texttt{lsize}=100$) for uniform rp16, uniform fp32, and \APIC  with target accuracy $\epsrel=10^{-4}$ or $\epsrel=10^{-8}$.
Values of this ratio that are less than 1 are an improvement. 
For $\texttt{lsize}=100$, \APIC achieves the lowest cost for eight of the problems, using $\epsrel=10^{-4}$ for seven of them
and $\epsrel=10^{-8}$ for one (Emilia\_923). For the remaining two problems (Hook\_1498 and boneS10), uniform rp16 is the cheapest, but \APIC with $\epsrel=10^{-4}$
has a similar cost, 
which shows that \APIC is automatically able to detect easy problems where reduced precision is sufficient and to match its cost.
In contrast, on the eight problems where \APIC does better than rp16, the cost difference can be significant; 
in fact, on some problems (ldoor and Emilia\_923), rp16 can be more expensive than fp64 (with cost up to $6.5\times$ larger), 
whereas \APIC with $\epsrel=10^{-4}$ always costs less than fp64 (relative cost of at most $0.82$). This shows that \APIC
only approximates when it is safe to do so and maintains the robustness of uniform high precision factorizations.

For $\texttt{lsize}=20$ the same trends are observed. This case is slightly more favourable for uniform rp16, which is the cheapest for four problems out of the nine that are successfully solved for this value of \texttt{lsize} (PCG fails for boneS10 regardless of the precision choice). This is once again consistent with previous observations that showed that for highly sparse factors, the effect of storing the entries in very reduced precision is less important. 
However, for two of these four problems (Flan\_1565 and audikw\_1),  \APIC with $\texttt{lsize}=100$ achieves a lower cost
than uniform rp16 with $\texttt{lsize}=20$; hence, across both values of \texttt{lsize}, \APIC is the cheapest variant
for seven of the ten problems.

As a side remark, we note that adaptive precision with $\epsrel=10^{-8}$ is not always cheaper in storage than fp32 (whose unit roundoff $u=2^{-24}$ is slightly larger than $\epsrel$) because entries that require more accuracy than that of fp32 are promoted to rp40.

\begin{table}[tbp]
  \centering
  \footnotesize
  \caption{\footnotesize Results for the ten test problems with the largest number of nonzeros in $A$ (given in brackets in millions), for $\texttt{lsize} = 20$ and $100$, for the uniform fp64, fp32, and rp16 factorizations, and the \APIC variants (shortened to AP in the table) with $\epsrel = 10^{-4}$ and $10^{-8}$. ``FP'' is the storage for the floating-point values in $\Lh$, and ``Total'' is the storage for $\Lh$ in sparse format, which adds 4-byte row indices and column pointers; both are in MB. ``\#its'' is the PCG iteration count, and ``R.Cost'' is the cost~\eqref{eqn:cost} normalized relative to that of uniform fp64 with $\texttt{lsize} = 100$. For each problem and \texttt{lsize}, the lowest R.Cost is in bold. A --- denotes PCG did not converge to the requested accuracy within the iteration limit.}
  \label{tab:memlim_largest}
  \begin{tabular}{llrrrrrrrr}
    \toprule
    & & \multicolumn{4}{c}{\texttt{lsize} $= 20$} & \multicolumn{4}{c}{\texttt{lsize} $= 100$} \\
    \cmidrule(lr){3-6}\cmidrule(lr){7-10}
    Problem & Method & FP & Total & \#its & R.Cost & FP & Total & \#its & R.Cost \\
    \midrule
    \multirow{5}{*}{\begin{tabular}[c]{@{}l@{}}Flan\_1565\\(59.5M)\end{tabular}} & rp16 & 173 & 525 & 1122 & \textbf{0.61} & 412 & 1242 & 606 & 0.60 \\
     & AP $10^{-4}$ & 209 & 561 & 1117 & 0.63 & 448 & 1277 & 570 & \textbf{0.57} \\
     & fp32 & 346 & 699 & 1118 & 0.72 & 824 & 1654 & 570 & 0.71 \\
     & AP $10^{-8}$ & 353 & 705 & 1118 & 0.73 & 821 & 1651 & 570 & 0.71 \\
     & fp64 & 693 & 1045 & 1119 & 0.97 & 1648 & 2477 & 570 & 1.00 \\
    \cmidrule(lr){1-10}
    \multirow{5}{*}{\begin{tabular}[c]{@{}l@{}}audikw\_1\\(39.3M)\end{tabular}} & rp16 & 111 & 336 & 432 & \textbf{0.66} & 255 & 768 & 289 & 0.78 \\
     & AP $10^{-4}$ & 137 & 362 & 524 & 0.84 & 281 & 794 & 213 & \textbf{0.59} \\
     & fp32 & 222 & 447 & 804 & 1.47 & 510 & 1023 & 209 & 0.71 \\
     & AP $10^{-8}$ & 226 & 452 & 788 & 1.45 & 505 & 1018 & 211 & 0.71 \\
     & fp64 & 444 & 669 & 790 & 1.92 & 1020 & 1533 & 209 & 1.00 \\
    \cmidrule(lr){1-10}
    \multirow{5}{*}{\begin{tabular}[c]{@{}l@{}}bone010\\(36.3M)\end{tabular}} & rp16 & 107 & 325 & 1272 & 0.74 & 257 & 776 & 688 & 0.73 \\
     & AP $10^{-4}$ & 132 & 349 & 1072 & \textbf{0.65} & 282 & 801 & 593 & \textbf{0.65} \\
     & fp32 & 214 & 431 & 992 & 0.69 & 515 & 1033 & 533 & 0.72 \\
     & AP $10^{-8}$ & 217 & 435 & 984 & 0.69 & 512 & 1031 & 524 & 0.70 \\
     & fp64 & 428 & 645 & 984 & 0.91 & 1030 & 1548 & 525 & 1.00 \\
    \cmidrule(lr){1-10}
    \multirow{5}{*}{\begin{tabular}[c]{@{}l@{}}Serena\\(33.0M)\end{tabular}} & rp16 & 116 & 353 & 97 & 0.49 & 328 & 990 & 63 & 0.69 \\
     & AP $10^{-4}$ & 144 & 381 & 84 & \textbf{0.44} & 356 & 1018 & 50 & \textbf{0.56} \\
     & fp32 & 232 & 469 & 85 & 0.52 & 656 & 1318 & 50 & 0.70 \\
     & AP $10^{-8}$ & 237 & 474 & 85 & 0.52 & 641 & 1302 & 50 & 0.69 \\
     & fp64 & 464 & 701 & 85 & 0.70 & 1313 & 1975 & 50 & 1.00 \\
    \cmidrule(lr){1-10}
    \multirow{5}{*}{\begin{tabular}[c]{@{}l@{}}Geo\_1438\\(32.3M)\end{tabular}} & rp16 & 116 & 355 & 94 & \textbf{0.54} & 336 & 1013 & 46 & 0.58 \\
     & AP $10^{-4}$ & 139 & 377 & 94 & 0.56 & 359 & 1036 & 43 & \textbf{0.55} \\
     & fp32 & 233 & 471 & 94 & 0.65 & 672 & 1349 & 43 & 0.70 \\
     & AP $10^{-8}$ & 234 & 473 & 94 & 0.65 & 641 & 1318 & 43 & 0.68 \\
     & fp64 & 466 & 704 & 94 & 0.88 & 1343 & 2021 & 43 & 1.00 \\
    \cmidrule(lr){1-10}
    \multirow{5}{*}{\begin{tabular}[c]{@{}l@{}}Hook\_1498\\(31.2M)\end{tabular}} & rp16 & 117 & 356 & 426 & \textbf{0.52} & 345 & 1041 & 195 & \textbf{0.54} \\
     & AP $10^{-4}$ & 141 & 380 & 427 & 0.55 & 370 & 1066 & 194 & 0.55 \\
     & fp32 & 233 & 472 & 427 & 0.64 & 690 & 1387 & 194 & 0.69 \\
     & AP $10^{-8}$ & 237 & 476 & 427 & 0.64 & 658 & 1354 & 194 & 0.68 \\
     & fp64 & 467 & 706 & 427 & 0.86 & 1381 & 2077 & 194 & 1.00 \\
    \cmidrule(lr){1-10}
    \multirow{5}{*}{\begin{tabular}[c]{@{}l@{}}boneS10\\(28.2M)\end{tabular}} & rp16 & 89 & 270 & --- & --- & 228 & 688 & 979 & \textbf{0.56} \\
     & AP $10^{-4}$ & 110 & 291 & --- & --- & 251 & 711 & 974 & 0.57 \\
     & fp32 & 177 & 358 & --- & --- & 457 & 916 & 975 & 0.70 \\
     & AP $10^{-8}$ & 181 & 362 & --- & --- & 440 & 900 & 975 & 0.69 \\
     & fp64 & 355 & 536 & --- & --- & 913 & 1373 & 974 & 1.00 \\
    \cmidrule(lr){1-10}
    \multirow{5}{*}{\begin{tabular}[c]{@{}l@{}}ldoor\\(23.7M)\end{tabular}} & rp16 & 82 & 248 & 274 & 2.51 & 227 & 684 & 253 & 4.94 \\
     & AP $10^{-4}$ & 100 & 267 & 156 & \textbf{1.50} & 245 & 702 & 30 & \textbf{0.60} \\
     & fp32 & 163 & 330 & 156 & 1.73 & 454 & 911 & 28 & 0.70 \\
     & AP $10^{-8}$ & 164 & 331 & 156 & 1.74 & 401 & 858 & 28 & 0.66 \\
     & fp64 & 326 & 493 & 156 & 2.34 & 907 & 1365 & 28 & 1.00 \\
    \cmidrule(lr){1-10}
    \multirow{5}{*}{\begin{tabular}[c]{@{}l@{}}Emilia\_923\\(21.0M)\end{tabular}} & rp16 & 75 & 229 & 339 & 3.22 & 216 & 652 & 314 & 6.57 \\
     & AP $10^{-4}$ & 82 & 236 & 65 & 0.63 & 223 & 659 & 39 & 0.82 \\
     & fp32 & 150 & 304 & 48 & 0.55 & 432 & 868 & 26 & 0.70 \\
     & AP $10^{-8}$ & 142 & 296 & 48 & \textbf{0.54} & 366 & 801 & 26 & \textbf{0.65} \\
     & fp64 & 301 & 455 & 48 & 0.75 & 864 & 1300 & 26 & 1.00 \\
    \cmidrule(lr){1-10}
    \multirow{5}{*}{\begin{tabular}[c]{@{}l@{}}inline\_1\\(18.7M)\end{tabular}} & rp16 & 55 & 166 & 1780 & 1.29 & 132 & 397 & 542 & 0.72 \\
     & AP $10^{-4}$ & 67 & 178 & 1680 & \textbf{1.26} & 143 & 409 & 422 & \textbf{0.57} \\
     & fp32 & 110 & 221 & 1646 & 1.43 & 263 & 529 & 422 & 0.71 \\
     & AP $10^{-8}$ & 112 & 223 & 1679 & 1.46 & 252 & 517 & 422 & 0.69 \\
     & fp64 & 219 & 331 & 1647 & 1.90 & 527 & 792 & 422 & 1.00 \\
    \bottomrule
  \end{tabular}
\end{table}

Finally, we consider the impact of problem size (in terms of the number of nonzero entries). In general, we  observe that the more aggressive factorizations (such as uniform rp16 or \APIC with a large $\epsrel$) perform better on the large test problems, 
outperforming uniform fp64 on a higher percentage of problems and with lower costs on average.
For example, with $\texttt{lsize}=20$, the median cost of uniform rp16 (normalized relative to uniform fp64) is $0.88$ for the problems in the largest third ($nz\ge 1.6\times10^{6}$) but is $1.9$ for the problems in the smallest third ($nz\le 8.9\times10^{4}$).
While the same trend is observed for \APIC, it is much less pronounced: with $\epsrel=10^{-4}$, the median cost increases from $0.64$ (for the larger problems) to $0.75$ (for the smaller problems).
Thus, while problem size also impacts the performance of \APIC, in our experiments it is much less sensitive to it than uniform rp16 and achieves significant gains even on the smaller problems of the set.

Overall, the performance profiles  (Figure~\ref{fig:memlim_perf_prof}) together
with the results on large problems (Table~\ref{tab:memlim_largest}) demonstrate that the use of adaptive precision in memory-limited IC factorizations can significantly reduce the size of the factors while maintaining fast convergence, leading to substantial reductions in the cost metric \eqref{eqn:cost}. As future work, this motivates developing an optimized high-performance implementation of \APIC (based, for example, on memory accessors) to assess to what extent these cost reductions  translate into effective runtime gains.

\section{Conclusion}
\label{sec:concl}

We have proposed adaptive precision incomplete LU and Cholesky factorizations
(Algorithms~\ref{alg.ILUT} and \ref{alg.IC})
in which the storage precision for each entry of the factors is chosen according to the target accuracy $\epsabs$
and its magnitude: entries of smaller magnitude can be stored in correspondingly lower precision.
Threshold-based factorizations can be directly incorporated into this analysis by considering dropped elements to be stored at the lowest ``precision'', and Theorems~\ref{thm.ILUT} and~\ref{thm.IC} show that this yields a backward error of order $\epsabs$ for the adaptive precision ILU and IC factorizations, respectively.
Our numerical examples illustrate that this can reduce the storage required for a given accuracy compared with uniform precision factors, including for matrices for which little dropping is possible.

We have also discussed several practical considerations.
Pivoting does not change the analysis or the precision selection criteria. Scaling can be accommodated either by applying the algorithm to the scaled matrix or by reverting the scaling in the precision selection criteria.
The latter can give larger storage savings, but we found the former to be more robust.
For triangular solves, performing the substitutions in high precision gives the standard backward error.
By using memory accessors, storage reduction can still be exploited with high precision solves if the solves are memory-bound.
Evaluating the solves in adaptive precision, by contrast, can be unstable in the presence of small pivots unless rook pivoting is used in the factorization, as our analysis and experiments demonstrate, 
or care is taken to convert to high precision entries in columns with small pivots.
We recommend keeping the solves in higher precision.

Finally, we have extended the approach to memory-limited factorizations by taking $\epsabs$ to be the magnitude of the largest dropped entry.
Experiments with a modified version of 
the IC software {\tt HSL\_MI28} on 153 SPD problems showed that, with a suitably chosen $\epsabs$, adaptive precision storage improves a metric designed to model cost (the product of the storage requirements and number of PCG iterations required to solve a system) compared with uniform fp64, fp32, and rp16 factors for the majority of the test set.
These experiments also showed that the appropriate choice of $\epsabs$ depends on the fill parameter \texttt{lsize}, which indirectly controls the largest dropped value.

Several directions remain for future work.
First, we have compared preconditioners using the cost model \eqref{eqn:cost};
a dedicated performance study of an efficient sparse implementation is needed to determine the extent to which reduced storage translates into shorter factorization and solve times on memory-bound architectures (Remark~\ref{rem.facto-time}).
Second, for memory-limited factorizations, the proposed appropriate $\epsabs$ is only known a posteriori.
Using a dynamic $\epsabs$ that changes during the factorization, for example as a running maximum of the dropped entries, would avoid the need for a prior estimate (Remark~\ref{rem.dynamic_epsa}).
Third, the computations in \APILU and \APIC could also be performed in adaptive precision (Remark~\ref{rem.AP-compute}).
Finally, it would be of interest to study the extension to other classes of incomplete factorization.

\section*{Acknowledgements}

This work was partially supported by projects
managed by the French National Research Agency (ANR):
France 2030 NumPEx Exa-MA (ANR-22-EXNU-0002),
MixHPC (ANR-23-CE46-0005),
and
FPT-4 (ANR-24-CE46-7572).

For the purpose of open access, the author has applied a Creative Commons Attribution (CC BY) licence to any Author Accepted Manuscript version arising.

\bibliographystyle{siamplain}
\bibliography{strings,refs,tmary}

\end{document}